\documentclass[12pt, reqno]{amsart}
\usepackage{amsmath,amssymb,amsthm}
\usepackage{amsmath, amssymb}
\usepackage{amsthm, amsfonts, mathrsfs}
\usepackage{mathptmx}
\usepackage{fullpage}
\usepackage{amsfonts,graphicx}
\usepackage{dutchcal}
\numberwithin{equation}{section}
\usepackage[colorlinks=true, pdfstartview=FitV, linkcolor=blue, citecolor=blue, urlcolor=blue]{hyperref}

\newtheorem {theorem}{Theorem}[section]
\newtheorem {lemma}[theorem]{{\bf Lemma}}
\newtheorem {corollary}[theorem]{{\bf Corollary}}

\theoremstyle{remark}
\newtheorem {remark}{{\bf Remark}}[section]

\theoremstyle{plain} \numberwithin {equation}{section}

\newcommand{\R}{{\mathbb R}}
\def\div{ \hbox{\rm div}\,  }
\newcommand\Z{{\mathbb{Z}}}

\def\nn{\nonumber}

\def\u{ \mathbf{u} }
\def\v{ \mathbf{v} }
\def\T{ \mathbb{T} }

\begin{document}
\title{ { Stability threshold for two-dimensional Boussinesq systems \\[1ex]near-Couette shear flow in a finite channel }}

\author{Tao Liang}

\address{ School of Mathematics,
	South China University of Technology,
	Guangzhou, 510640, China}

\email{taolmath@163.com}

\author{Yongsheng Li}
\address{ School of Mathematics,
	South China University of Technology,
	Guangzhou, 510640, China}

\email{yshli@scut.edu.cn}

\author{Xiaoping Zhai}

\address{ Department of Mathematics, Guangdong University of Technology,
	Guangzhou, 510520, China}

\email{pingxiaozhai@163.com (Corresponding author)}




\begin{abstract}
In this paper, we investigate the stability threshold problem of the two-dimensional Navier-Stokes Boussinesq(NSB) equations in a finite channel $ \T \times [-1,1]$, focusing on the stability around the near Couette shear flow $ (U(y), 0)$, assuming the Navier slip boundary conditions are satisfied. In particular, when the initial data for the vorticity resides in an anisotropic Sobolev space of size $ O(\min \{ \mu^{\frac{1}{2}}, \nu^{\frac{1}{2}}\})$, and the initial perturbation of the temperature resides in an anisotropic Sobolev space of size $ O(\min \{ \mu, \nu\})$, we derive the nonlinear enhanced dissipation effect and the inviscid damping effect for the NSB system.
\end{abstract}

\maketitle

\section{Introduction and the main results}

In this paper, we consider the two-dimensional Navier-Stokes Boussinesq system with Navier slip boundary conditions over a finite channel $(x,y) \in \T \times [-1, 1]$, given as follows:
\begin{eqnarray}\label{quanhs}
\left\{\begin{aligned}
&\partial_t \v+ \v \cdot\nabla \v - \mu (\sigma \partial_x^2 + \partial_y^2) \v  + \nabla P  = \rho g e_2  ,\\
& \partial_t \rho + \v \cdot\nabla\rho - \nu (\sigma \partial_x^2 + \partial_y^2) \rho  = 0,\\
& \div \v = \partial_x v_1 + \partial_y v_2 = 0,\\
& v_2(t, x, \pm 1) = 0, \quad \partial_y v_1(t, x, \pm 1) = 1, \quad \rho(t, x, \pm 1) = c_0,\\
& \v(0, x,y) = \v_{in}(x,y),\qquad \rho(0,x,y) = \rho_{in}(x, y).
\end{aligned}\right.
\end{eqnarray}
The variables $\v$, $\rho$, and $P$ denote the velocity, the temperature, and the pressure, respectively. Here, $ \mu$, $ \nu$ and $g$ represent the viscosity coefficient, thermal diffusivity coefficient, and the  gravitational constant, respectively, while $e_2 = (0,1)^\top$ denotes the unit vector in the vertical direction. In System \eqref{quanhs}, the boundary conditions imply that the fluid is allowed to slip freely near the boundary, while  $\rho(t, x, \pm 1) = c_0$ represents a fixed temperature boundary condition. Combining these two velocity boundary conditions with the temperature boundary condition, the system likely models an interaction between low-friction flow and thermal convection. In order to simplify the discussion without loss of generality, we normalize the constants such that $ g = c_0 =1$.

In this paper, we consider the stability problem in the vicinity of the following steady-state solution
\begin{align*}
	\v_s = ( U(t,y), 0)^{\mathrm{T}}, \quad  \rho_s = 1, \quad P_s = y + c.
\end{align*}
Here, note that $$ U(t,y) = y + \partial_y \int_{-1}^{1} G(y, y') W(t, y') \ \mathrm{d}y'$$ is the shear flow determined by the Biot-Savart law, where
\begin{align*}
	 G\left(y, y^{\prime}\right)= \begin{cases}\frac{-\left(y^{\prime}-1\right)}{2}(y+1), & \text { for } y\le y^{\prime}; \\ \frac{-\left(y^{\prime}+1\right)}{2}(y-1), & \text { for } y\ge y^{\prime};\end{cases}
\end{align*}
 is the Green's function determined by $ \partial_y^2 G(y, y') = \delta(y-y')$ with homogeneous Dirichlet boundary conditions (for $y\in [-1,1]$); and $ W(t,y) = e^{\mu\partial_y^2 t}W_{in}(y) $  (for $y\in[-1,1]$) is the solution to the heat equation
 \begin{equation*}
 \left\{\begin{array}{l}
 \partial_t W(t,y) - \mu \partial_y^2 W(t, y) = 0, \\
 W|_{t=0} = W_{in}(y), \qquad W(t,  \pm 1)=0.
 \end{array}\right.
 \end{equation*}

Let $ \u = \v-\v_s, \ \theta = \rho - \rho_s, \ p = P - P_s$,  then system \eqref{quanhs} can  be rewritten as the equation for $ (\u, \theta, p)$ as follows:
\begin{eqnarray}\label{rewrite}
\left\{\begin{aligned}
&\partial_t \u+ U\partial_x \u + \binom{u_2 \partial_y U}{ 0} + \u \cdot\nabla \u - \mu (\sigma \partial_x^2 + \partial_y^2) \u  + \nabla p  =  \binom{0}{ \theta},\\
& \partial_t \theta + U\partial_x \theta + \u \cdot\nabla\theta - \nu (\sigma \partial_x^2 + \partial_y^2) \theta  = 0,\\
& \div \u =  0,\\
& u_2(t, x, \pm 1) = 0, \quad \partial_y u_1(t, x, \pm 1) = 0, \quad \theta(t, x, \pm 1) = 0.
\end{aligned}\right.
\end{eqnarray}
By introducing the vorticity $ \omega = \nabla \times \u = \partial_y u_1 - \partial_x u_2$, the above system \eqref{rewrite} can be reformulated as follows
\begin{eqnarray}\label{rewrite1}
\left\{\begin{aligned}
&\partial_t \omega+ U\partial_x \omega - U''\partial_x \psi + \u \cdot\nabla \omega - \mu (\sigma \partial_x^2 + \partial_y^2) \omega    =  -\partial_x \theta,\\
& \partial_t \theta + U\partial_x \theta + \u \cdot\nabla\theta - \nu (\sigma \partial_x^2 + \partial_y^2) \theta  = 0,\\
& \omega(0, x, y) = \omega_{in}(x,y), \quad \theta(0, x, y) = \theta_{in}(x,y), \\
& \omega(t, x, \pm 1) = 0, \quad  \theta(t, x, \pm 1) = 0, \\
& \u = \nabla^{\perp} \psi = (\partial_y \psi,-\partial_x \psi ), \quad \Delta \psi =\omega.
\end{aligned}\right.
\end{eqnarray}
We investigate the following stability threshold problem.
\vskip .1in
 Given norms $\|\cdot\|_{Y_1}$  and $\|\cdot\|_{Y_2}$, find $\alpha=\alpha\left(Y_1, Y_2\right) $ and  $\beta=\beta\left(Y_1, Y_2\right)$  such that
\begin{itemize}
  \item  $\left\|\omega_{\mathrm{in}}\right\|_{Y_1} \leq \mu^\alpha$ \hbox{ and }$\left\|\theta_{\mathrm{in}}\right\|_{Y_2} \leq \nu^\beta \Rightarrow$ \hbox{stability };\\
  \item $\left\|\omega_{\mathrm{in}}\right\|_{Y_1} \gg \mu^\alpha$ \text { or }$\left\|\theta_{\mathrm{in}}\right\|_{Y_2} \gg \nu^\beta \Rightarrow$ \text { instability }.
\end{itemize}

\subsection {Background and recent studies}
The Navier-Stokes Boussinesq system serves as a fundamental model for buoyancy-driven flows and finds extensive applications in atmospheric sciences, oceanic circulation, mantle convection, and engineering heat transfer \cite{Gill,Majda}. Its core assumption is the Boussinesq approximation, which retains density variations only in the buoyancy term (such as temperature or salinity fluctuations) while neglecting them in the inertial term. This model effectively captures key phenomena such as Rayleigh-B\'enard convection and stratified shear flows \cite{Constantin,Getling}. Moreover, the vortex dynamics in 2D flows differ significantly from those in 3D: energy cascades inversely towards larger scales, while the coupling between buoyancy effects and viscous dissipation influences turbulence structures \cite{Kraichnan}. Additionally, the 2D NSB system retain some crucial features of the 3D Euler and Navier-Stokes equations, such as the vortex stretching mechanism, which plays a significant role in the nonlinear dynamics and energy transfer processes \cite{Majada2}.

Before stating our results, we first review previous studies on the stability of perturbations around steady states for System \eqref{quanhs}. The stability analysis is generally classified into two main categories based on the domain: one for $\T \times \R$ and the other for $\T \times [-1,1]$.

$(\mathrm{I})$ The domain $\T \times \R$.
\begin{itemize}
  \item When considering perturbations around the shear flow $ \v_s = (U(y),0)$ and in the absence of hydrostatic equilibrium $ \rho_s = 1 $, Tao and Wu \cite{Tao2017JDE} first established the linear stability of System \eqref{rewrite1} in the upper half-space $ \T \times \R^+$. Notably, they demonstrated that enhanced dissipation effects can be obtained with dissipation present only in the vertical direction ($ \sigma= 0$).
 Subsequently, Deng {\it et al.} \cite{Deng2021JFA} used Fourier multipliers to establish the nonlinear stability of System \eqref{rewrite1} under the same vertical dissipation assumption ($ \sigma = 0 $ and $ \mu = \nu$), provided that the initial perturbation satisfy $$ \| \omega_{in}\|_{H^b} \lesssim \nu^{\frac{2}{3}}, \quad \text{and} \quad \| \theta_{in}\|_{H^b} + \|\langle \partial_x\rangle^{\frac{1}{3}} \theta_{in}\|_{H^b} \lesssim \nu,$$ with $ b >\frac{4}{3}$.
  Additionally, they established the nonlinear stability of System \eqref{rewrite1} under full dissipation when the initial perturbation satisfy $$\| \omega_{in}\|_{H^b} \lesssim \nu^{\frac{1}{2}}, \quad \text{and} \quad \| \theta_{in}\|_{H^b} + \nu^{\frac{1}{6}}\|\langle \partial_x\rangle^{\frac{1}{3}} \theta_{in}\|_{H^b} \lesssim \nu,$$ with $ b > 1$, thereby improving the stability threshold.
   Their results,  mathematically reveal the mechanism of inviscid damping and enhanced dissipation. Recently, Zhang and Zi \cite{Zhang2023JMPA} further refined the stability threshold in \cite{Deng2021JFA} under the assumption that the initial perturbation satisfy $$ \| \omega_{in}\|_{H^b}   + \nu^{-\frac{1}{2}}\| \theta_{in}\|_{H^b} + \nu^{-\frac{1}{3}}\|\langle \partial_x\rangle^{\frac{1}{3}} \theta_{in}\|_{H^b} \lesssim \nu^{\frac{1}{3}},$$ with $ b > 7$.
Meanwhile, Niu and Zhao \cite{Niu2024ARxiv} combined the quasi-linearization method from \cite{Chenqi2020ARMA,Chenqi2024AMS} with the time-dependent elliptic operator $ \Lambda_t^2 = 1-\partial_x^2 - (\partial_y + t\partial_x)^2$ constructed in \cite{Deng2021JFA,Wei2023Tunisi} and employed Fourier analysis to further improve the stability results of \cite{Zhang2023JMPA}. Notably, they raised the threshold for the initial temperature perturbation to $ \nu^{\frac{2}{3}}$. To the best of our knowledge, this is the optimal stability threshold currently available when $ \mu = \nu$.
 \vskip .1in
  \item When considering perturbations around shear flow  with hydrostatic equilibrium, i.e., $ \v_s = ( U(y), 0)$ and $ \rho_s = \rho_s(y)$, Yang and Lin \cite{Lin2018JMFM} first established the linear stability of system \eqref{quanhs} in the inviscid case for the Couette flow $ \v_s = ( y, 0)$ in the inviscid case, where $ \rho_s (y) $ represents an exponentially stratified density profile. Subsequently, Bianchini {\it et al.} \cite{Dolce2021India} extended this result to the linear stability near Couette flow $ \v_s = ( U(y), 0)$($ \| U(y) - y\|_{H^s} \ll 1$). Later, under appropriate constraints on the Richardson number, Bedrossian {\it et al.} \cite{Dolce2023CPAM} demonstrated that, given an initial perturbation of Gevrey class with magnitude $ \varepsilon$, the system exhibits shear-buoyancy instability.  Moreover, Masmoudi {\it et al.} \cite{Zhao2022ARMA} proved the stability of initial data in Gevrey$-\frac{1}{s} (\frac{1}{3} < s \le 1) $ in the absence of thermal diffusion. As for the system with full dissipation, Zillinger obtained the nonlinear stability and enhanced dissipation in \cite{Zillinger1}. The nonlinear results was extended by Zillinger himself in \cite{Zillinger2} to possible oscillating temperature profiles for system with only vertical viscous and thermal dissipation.
  Additionally, Zhai and Zhao \cite{Zhao2023SIAM} recently investigated the nonlinear asymptotic stability of the Couette flow under certain constraints on the Richardson number, further extending the understanding of stability in this setting.
\end{itemize}

$(\mathrm{II})$ The domain $\T \times [-1,1]$.
 When considering perturbations around the shear flow $ \v_s = (U(y),0)$ and in the absence of hydrostatic equilibrium $ \rho_s = 1 $,  Masmoudi {\it et al.} \cite{Masmoudi2023JFA} established the nonlinear stability of system \eqref{rewrite1} for the Couette shear flow by constructing suitable weighted energy norms and leveraging the resolvent estimates of the linearized operator from \cite{Chenqi2020ARMA}, along with space-time estimates for the linearized Navier-Stokes equations. This result further advances the understanding of the long-term dynamical behavior of the Boussinesq system in a finite domain, particularly in relation to stability mechanisms under small initial perturbations.

\subsection{Main results}
First, we analyze the linearized stability of System \eqref{rewrite1}.
To make the statement precise, we define, for $(x,y) \in \T \times [-1, 1]$,
 $$ f_k(t,y) = \frac{1}{2\pi} \int_{-\pi}^{\pi} f(t,x,y)e^{-ikx} \mathrm{d}x.$$
Taking a Fourier expansion in the periodic direction $x$, we obtain  the modes of $ (\omega, \theta)$ for each wave number $k$ as follows:
\begin{eqnarray}\label{linear_eq}
\left\{\begin{aligned}
&\partial_t \omega_k+ ik U \omega_k - ik U'' \psi_k - \mu (-\sigma k^2 + \partial_y^2) \omega_k    =  -i k \theta_k,\\
& \partial_t \theta_k + i k U \theta_k  - \nu (-\sigma k^2 + \partial_y^2) \theta_k  = 0,\\
&    \Delta_{k} \psi_k =(-k^2 + \partial_y^2)\psi_k=\omega_k,\\
& \psi_k (t, \pm 1) = \omega_k(t, \pm 1) = 0, \quad  \theta_k(t, \pm 1) = 0.
\end{aligned}\right.
\end{eqnarray}

Now, the   first main result of this paper is stated as follows.
\begin{theorem}(Linearized stability)
\label{theorem1.1}
Set $\sigma = 0$ and let
$ (\omega, \theta)$ denote the solution to the linearized equations associated with System \eqref{rewrite1}. Under the initial boundary conditions $ W_{in}|_{y = \pm 1} = \omega_{in}|_{y = \pm 1} = 0$, there exists a small positive constant $\delta_0 \in (0,\frac{1}{2})$ such that
 if the initial data satisfies
\begin{align}\label{W0}
& \|W_{in}(\cdot) \|_{H^4_y}  \leq \delta_0,
\end{align}
then for all integers $k \neq 0$ and any time $t>0$, the following estimates hold:
\begin{align}\label{linear_ta}
	\| \theta_k\|^2_{L^2} + \nu^{\frac{2}{3}} |k|^{-\frac{2}{3}} \| \partial_y \theta_k\|^2_{L^2} \le C_0 \Phi_{\theta_{in}} e^{-c_0 \nu^{\frac{1}{3}|k|^{\frac{2}{3}} t}},
\end{align}
and
\begin{align}\label{linear_w}
	\| \omega_k\|^2_{L^2} + \mu^{\frac{2}{3}} |k|^{-\frac{2}{3}} \| \partial_y \omega_k\|^2_{L^2} \le C_0 \Phi_{\omega_{in}} e^{-c_0 \mu^{\frac{1}{3}|k|^{\frac{2}{3}} t}}+ C_0 \min\{ \mu, \nu  \}^{-\frac{2}{3}} |k|^{\frac{2}{3}}  \Phi_{\theta_{in}} e^{-c_0 \min\{\mu, \nu\}^{\frac{1}{3}|k|^{\frac{2}{3}} t}},
\end{align}
where $ C_0$ and $ c_0$ are two constants independent of $k$ and $t$. Here $ \Phi_{\omega_{in}}$ and $ \Phi_{\theta_{in}}$  are defined as
\begin{align*}
	&\Phi_{\omega_{in}} = \| \omega_{in}(k, y)\|^2_{L^2} + \mu^{\frac{2}{3}} |k|^{-\frac{2}{3}} \| \partial_y \omega_{in}(k, y)\|^2_{L^2},\quad
	& \Phi_{\theta_{in}} = \| \theta_{in}(k, y)\|^2_{L^2} + \nu^{\frac{2}{3}} |k|^{-\frac{2}{3}} \| \partial_y \theta_{in}(k, y)\|^2_{L^2}.
\end{align*}
\end{theorem}
\vskip .1in
\begin{remark}
Theorem \ref{theorem1.1} generalizes the seminal result of Tao and Wu \cite{Tao2017JDE} originally established on the domain $\T \times \R^+$ to the present framework. This extension not only broadens the applicability of their methodology but also addresses key technical challenges inherent to the current geometric and functional setting.
\end{remark}
\begin{remark}
We emphasize that the linear stability of System \eqref{rewrite1} requires only vertical dissipation in our framework. This crucial observation highlights the minimal dissipation requirements for stabilizing perturbations. Furthermore, our methodology naturally extends to systems with full dissipation (i.e., both vertical and horizontal dissipation), as the proof for the vertical dissipation case remains fully applicable without additional technical barriers.
\end{remark}

\vskip .1in

The main focus of this paper is actually the nonlinear stability.
The  second main result of this paper is stated as follows.
\begin{theorem}\label{thm1.2}(Nonlinear stability)
	Let $\sigma = 1$, and consider the initial boundary condition $$ W_{in}|_{y = \pm 1} = \omega_{in}|_{y = \pm 1} = \theta_{in}|_{y = \pm 1} = 0,$$ there exist  small positive constants $\delta_0, \varepsilon_0, \varepsilon_1$ such that, if
\begin{align*}
& \|W_{in}(\cdot) \|_{H^4_y}  \leq \delta_0,
\end{align*}
and
	\begin{align*}
		& \sum_{0\le j\le 1 } \big\| (\mu^{\frac{1}{3} } \partial_{y})^j \langle \partial_x\rangle^{m-\frac{j}{3}} \omega_{in} \big\|_{L^2} \le \varepsilon_0 \min\big\{ \mu^{\frac{1}{2} }, \nu^{\frac{1}{2} }\big\},\\
&\sum_{0\le j\le 1 } \big\| (\nu^{\frac{1}{3} } \partial_{y})^j \langle \partial_x\rangle^{m-\frac{j}{3}} \theta_{in} \big\|_{L^2} \le \varepsilon_1 \min\big\{ \mu, \nu \big\}.
	\end{align*}
Then, for any sufficiently small positive constant
$\delta_1$ (independent of
$\nu$ and $\mu$) and for all parameters
$ \mu, \nu \in (0,1)$, System \eqref{rewrite1} admits a globally well-posed solution
$ (\omega, \theta)$ satisfying the following uniform stability estimates:
	\begin{align*}
		\sum_{0\le j\le 1 } \big\| (\mu^{\frac{1}{3} } \partial_{y})^j  \omega_{0}  \big\|_{L^2} \le & C \varepsilon_0  \min\big\{ \mu^{-\frac{1}{2} }, \nu^{-\frac{1}{2} }\big\} e^{-\delta_1 \min\big\{ \mu, \nu\big\} t};\\
			\sum_{0\le j\le 1 } \big\| (\nu^{\frac{1}{3} } \partial_{y})^j  \theta_{0} \big\|_{L^2} \le & C \varepsilon_1  \min\big\{ \mu^{-1 }, \nu^{-1 }\big\} e^{-\delta_1 \min\big\{ \mu, \nu\big\} t};\\
		\sum_{0\le j\le 1 } \big\| (\mu^{\frac{1}{3} } \partial_{y})^j \langle \partial_x\rangle^{m-\frac{j}{3}} \omega_{\neq} \big\|_{L^2} \le & C \varepsilon_0 \min\big\{ \mu^{-\frac{1}{2} }, \nu^{-\frac{1}{2} }\big\} e^{-\delta_1 \min\big\{ \mu^{\frac{1}{3}}, \nu^{\frac{1}{3}} \big\} t};\\
		\sum_{0\le j\le 1 } \big\| (\nu^{\frac{1}{3} } \partial_{y})^j \langle \partial_x\rangle^{m-\frac{j}{3}} \theta_{\neq} \big\|_{L^2} \le & C\varepsilon_1 \min\big\{ \mu^{-1 }, \nu^{-1 }\big\} e^{-\delta_1 \min\big\{ \mu^{\frac{1}{3}}, \nu^{\frac{1}{3}} \big\} t},
	\end{align*}
for any  $t \in [0, +\infty)$.
\end{theorem}

\begin{remark}
In contrast to the hydrodynamic stability framework developed by Masmoudi {\it et al.} \cite{Masmoudi2023JFA}, our analysis achieves two critical advancements:
\begin{itemize}
  \item {\bf  Stability Threshold:} We establish a stability criterion for shear flows near the Couette profile, requiring the allowable initial perturbation $ \theta_{in}$ to be raised from $\frac{11}{12} $  to the  threshold of $1$. 
  \item {\bf Mechanism Decoupling:} We rigorously establish the coexistence of nonlinear inviscid damping and enhanced dissipation, notably accommodating distinct dissipation coefficients: $\mu$ (viscosity) and $\nu$ (thermal diffusivity).
\end{itemize}
\end{remark}
\vskip .1in
\subsection*{Notations}
Throughout this paper, $ C_0 > 0 $ denotes a generic constant that is independent of the relevant quantities. For brevity, we use the notation $ f \lesssim g $ to indicate that $  f \leq C_0 g$ for some constant $C_0 $. Let $ A$ and $B $ be two operators. We denote by $ [A, B] = AB - BA$ the commutator between $ A$ and $B $. Additionally, we use $ \langle f, g \rangle$ to represent the $ L^2_y([-1, 1])$ inner product of $ f$  and $g $.

\section{Linear Stability}\label{sec2}

\subsection{Preliminaries}
First, we construct two coercive energy functional (the constant $ c_{\alpha} \in (0,1)$ remains to be determined).
\begin{align}\label{Ekw}
E_{\omega, k}[\omega_k] = 128\| \omega_k\|^2_{L^2} + 4 \mu^{\frac{2}{3}} |k|^{-\frac{2}{3}} \|\partial_y \omega_k \|^2_{L^2} +  \mu^{\frac{1}{3}} |k|^{-\frac{4}{3}}
\operatorname{Re} \langle i k \omega_k, \partial_y \omega_k \rangle + c_{\alpha} \operatorname{Re} \langle  \omega_k, \mathfrak{J}_k [\omega_k] \rangle
\end{align}
and
\begin{align}\label{Ekta}
E_{\theta, k}[\theta_k] = 16\| \theta_k\|^2_{L^2} +  \nu^{\frac{2}{3}} |k|^{-\frac{2}{3}} \|\partial_y \theta_k \|^2_{L^2} +  \nu^{\frac{1}{3}} |k|^{-\frac{4}{3}}
\operatorname{Re} \langle i k \theta_k, \partial_y \theta_k \rangle
\end{align}
to examine the asymptotic behavior of $(\omega_k, \theta_k) $ for each mode $k$.The singular integral operator $ \mathfrak{J}_k$  in $ E_{\omega, k}[\omega_k]$ is defined as follows
\begin{align*}
\mathfrak{J}_k[f](y) :=|k| \mathrm{p} . \mathrm{v} . \frac{k}{|k|} \int_{-1}^1 \frac{1}{2 i(y-y^{\prime})} G_k(y, y^{\prime}) f(y^{\prime}) \mathrm{d} y^{\prime},
\end{align*}
where $G_k(y, y')$ denotes the Green's function determined by $ \Delta_{k} G(y, y') = \delta(y-y')$ with homogeneous Dirichlet boundary conditions. The exact expression for $ G_k(y, y')$ is given by
\begin{align*}
G_k\left(y, y^{\prime}\right)=-\frac{1}{k \sinh (2 k)} \begin{cases}\sinh \left(k\left(1-y^{\prime}\right)\right) \sinh (k(1+y)), & y \leq y^{\prime} ; \\ \sinh (k(1-y)) \sinh \left(k\left(1+y^{\prime}\right)\right), & y \geq y^{\prime}.\end{cases}
\end{align*}

Let us briefly recall the properties of the singular integral operator $ \mathfrak{J}_k$ introduced in \cite{WangCMP2025}( with $k \neq 0$ in the following three lemmas). First, we demonstrate that $ \mathfrak{J}_k$ is a bounded linear operator from $L^2$ to $L^2$, with the detailed proof provided in \cite{WangCMP2025}.

\begin{lemma}
\label{le2.1}
The singular integral operator $\mathfrak{J}_k$ is a bounded linear operator on $L^2 \rightarrow L^2$ and moreover
\begin{align*}
	\left\|\mathfrak{J}_k\right\|_{L^2 \rightarrow L^2} \lesssim 1.
\end{align*}
\end{lemma}
The next lemma captures the estimate for the commutator $ [\partial_y, \mathfrak{J}_k]$.
\begin{lemma}\label{le2.2}
	For the commutator $[\partial_y, \mathfrak{J}_k]$, there holds
	\begin{align*}
		\| [\partial_y, \mathfrak{J}_k] \|_{L^2 \rightarrow L^2} \lesssim |k|.
	\end{align*}
\end{lemma}

Finally, we describe the conjugate symmetry and self-adjointness of the operator $ \mathfrak{J}_k$.
\begin{lemma}\label{le2.3}
	For all $f, g \in L^2$ there holds
	\begin{align*}
		\overline{\mathfrak{J}_k[f]}=-\mathfrak{J}_k[\bar{f}],
	\end{align*}
	and
	\begin{align*}
	\int_{-1}^1 \bar{f} \mathfrak{J}_k[g] \mathrm{d} y=-\int_{-1}^1 \mathfrak{J}_k[\bar{f}] g \mathrm{d} y = \int_{-1}^1 \overline{\mathfrak{J}_k[f]} g \mathrm{d} y .
	\end{align*}
	In particular, we have $\mathfrak{J}_k=\mathfrak{J}_k^*$.
\end{lemma}

\vskip .1in

\subsection{Enhanced Dissipation Effect For $ \theta_k$}
In this subsection, we establish the enhanced dissipation effect for $ \theta_k$ using the coercive energy functional constructed from \eqref{Ekta}. Specifically, we obtain the following lemma.
\begin{lemma}\label{le_ta1}
	Under the conditions of Theorem \ref{theorem1.1}, we have the following:
	\begin{align}\label{1ta}
	\frac{d}{dt} E_{\theta, k} + \frac{1}{4} \nu^{\frac{1}{3}} |k|^{\frac{2}{3}} \|  \theta_k\|^2_{L^2} + \nu \| \partial_y \theta_k \|^2_{L^2} +  2\nu^{\frac{5}{3}} |k|^{-\frac{2}{3}} \|\partial_{yy} \theta_k \|^2_{L^2} \le 0.
	\end{align}
	\begin{proof}
		For the first two fundamental energy functionals in $ E_{\theta,k}[\theta_k]$, we can directly obtain
		\begin{align}\label{1ta_1}
		\frac{1}{2} \frac{d}{dt} \| \theta_k\|^2_{L^2} + \nu \| \partial_y \theta_k\|^2_{L^2} = 0,
		\end{align}
		and
		\begin{align}\label{1ta_two}
			\frac{1}{2} \frac{d}{dt} \nu^{\frac{2}{3}} |k|^{-\frac{2}{3}} \|\partial_y \theta_k \|^2_{L^2} + \nu^{\frac{5}{3}} |k|^{-\frac{2}{3}} \|\partial_{yy} \theta_k \|^2_{L^2} = -\nu^{\frac{2}{3}} |k|^{-\frac{2}{3}} \operatorname{Re} \langle ikU' \theta_k, \partial_y \theta_k \rangle.
		\end{align}
		According to the definition of $ U$ and the Green's function $ G$, there holds
		\begin{align*}
			U'(t,y) = 1+ \int_{-1}^{1} \partial_{yy} G(y,y') W(t,y') \mathrm{d} y' = 1+ W(t,y).
		\end{align*}
		It follows from the heat equation satisfied by $ W(t,y)$ and the homogeneous Dirichlet boundary conditions that we can readily obtain
		\begin{align}\label{W_estimate}
			\| W(t,y)\|_{H^4} \le \| W_{in}\|^2_{H^4} \le \delta_0.
		\end{align}
		For the local term in \eqref{1ta_two}, we obtain the following two estimates using Young's inequality
		\begin{align}\label{1ta_one}
		\big| \nu^{\frac{2}{3}} |k|^{-\frac{2}{3}} \operatorname{Re} \langle ikU' \theta_k, \partial_y \theta_k \rangle \big|  &\le \nu^{\frac{2}{3}} |k|^{\frac{1}{3}} \|U' \|_{L^\infty} \| \theta_k\|_{L^2} \|\partial_y \theta_k \|_{L^2}\nn\\
		&\le 2\nu^{\frac{2}{3}} |k|^{\frac{1}{3}} \| \theta_k\|_{L^2} \|\partial_y \theta_k \|_{L^2}\nn\\
		& \le 4\nu \|\partial_y \theta_k \|^2_{L^2} + \frac{1}{4}\nu^{\frac{1}{3}} |k|^{\frac{2}{3}} \| \theta_k\|^2_{L^2},
		\end{align}
and
		\begin{align}\label{1ta_2}
			\frac{1}{2} \frac{d}{dt} \nu^{\frac{2}{3}} |k|^{-\frac{2}{3}} \|\partial_y \theta_k \|^2_{L^2} + \nu^{\frac{5}{3}} |k|^{-\frac{2}{3}} \|\partial_{yy} \theta_k \|^2_{L^2} \le 8\nu \|\partial_y \theta_k \|^2_{L^2} + \frac{1}{8}\nu^{\frac{1}{3}} |k|^{\frac{2}{3}} \| \theta_k\|^2_{L^2}.
		\end{align}
		Next, we consider the following cross terms to generate the enhanced dissipation effect, we have
		\begin{align}\label{1ta_four}
			\frac{d}{dt} \operatorname{Re} \langle ik\theta_k, \partial_y \theta_k \rangle  &= \operatorname{Re} \langle ik(-ik U \theta_k + \nu \partial_{yy} \theta_k), \partial_y \theta_k \rangle + \operatorname{Re} \langle ik \theta_k, \partial_y (-ik U \theta_k + \nu \partial_{yy} \theta_k) \rangle\nn\\
			& = -\frac{|k|^2}{2} \| \sqrt{U'} \theta_k \|^2_{L^2} + \operatorname{Re} \langle ik\nu \partial_{yy} \theta_k, \partial_y \theta_k \rangle + \operatorname{Re} \int_{-1}^{1} ik \theta_k \cdot \overline{ -ik U' \theta_k}  \mathrm{d} y\nn\\
			& \quad  + \operatorname{Re} \int_{-1}^{1} ik \theta_k \cdot \overline{ -ik U \partial_y \theta_k}  \mathrm{d} y + \operatorname{Re} \int_{-1}^{1} ik \theta_k \cdot \overline{ \nu \partial_{yy} \partial_y \theta_k}  \mathrm{d} y\nn\\
			&  = -|k|^2 \| \sqrt{U'} \theta_k \|^2_{L^2} - \operatorname{Im} k\langle \nu \partial_{yy} \theta_k, \partial_y \theta_k \rangle - \operatorname{Im} k\langle \nu \partial_{yy} \partial_y\theta_k,  \theta_k \rangle.
		\end{align}
		By applying integration by parts ($ \theta_k|_{y = \pm 1} = 0$) and H\"older's inequality, we obtain
		\begin{align}\label{1ta_five}
		\nu^{\frac{1}{3}} |k|^{-\frac{4}{3}}\big| \operatorname{Im} k\langle \nu \partial_{yy} \theta_k, \partial_y \theta_k \rangle + \operatorname{Im} k\langle \nu \partial_{yy} \partial_y\theta_k,  \theta_k \rangle \big| &\le 2 \nu^{\frac{4}{3}} |k|^{-\frac{1}{3}}	\| \partial_y\theta_k\|_{L^2} \| \partial_{yy} \theta_k\|_{L^2}\nn\\
		&\le \nu \| \partial_y \theta_k\|^2_{L^2} + \nu^{\frac{5}{3}} |k|^{-\frac{2}{3}} \|\partial_{yy} \theta_k \|^2_{L^2}.
		\end{align}
		Collecting \eqref{1ta_four} and \eqref{1ta_five} yields
		\begin{align}\label{1ta_3}
		\nu^{\frac{1}{3}} |k|^{-\frac{4}{3}} \frac{d}{dt} \operatorname{Re} \langle ik\theta_k, \partial_y \theta_k \rangle + \nu^{\frac{1}{3}} |k|^{\frac{2}{3}} \| \sqrt{U'} \theta_k\|^2_{L^2} \le \nu \| \partial_y \theta_k \|^2_{L^2} +  \nu^{\frac{5}{3}} |k|^{-\frac{2}{3}} \|\partial_{yy} \theta_k \|^2_{L^2}.
		\end{align}
		Finally, multiplying \eqref{1ta_1} by 32, \eqref{1ta_2} by 2, and adding \eqref{1ta_3} give \eqref{1ta}.
	\end{proof}
\end{lemma}

\subsection{Enhanced Dissipation Effect For $ \omega_k$}
In this subsection, we establish the enhanced dissipation effect for $ \omega_k$ using the coercive energy functional constructed from \eqref{Ekw}.  The main difficulty arises from the presence of the local term $ -iU''\psi_k$ in the vorticity equation. To handle this, we introduce the singular integral operator $ \mathfrak{J}_k$ from \cite{WangCMP2025}. Additionally, for convenience, we define the energy of some dissipative terms as follows:
\begin{align*}
	&\operatorname{Dis}_{\omega, 1} = \mu \| \partial_y \omega_k\|^2_{L^2}, \quad \operatorname{Dis}_{\omega, 2} = \mu^{\frac{5}{3}} |k|^{-\frac{2}{3}} \| \partial_{yy} \omega_k\|^2_{L^2},\nn\\
	&\operatorname{Dis}_{\omega, 3} = \mu^{\frac{1}{3}} |k|^{\frac{2}{3}} \|  \omega_k\|^2_{L^2},\quad \operatorname{Dis}_{\omega, 4} = |k|^2 \| \nabla_k \psi_k\|^2_{L^2} \quad (\text{where} \nabla_k = (ik, \partial_{y})^\top ).
\end{align*}
Next, we proceed to estimate the energy functional defined in \eqref{Ekw} term by term.
\begin{lemma}\label{le_w1}
	Under the conditions of Theorem \ref{theorem1.1}, we have
	\begin{align}\label{1w}
		\frac{1}{2} \frac{d}{dt} \| \omega_k\|^2_{L^2} + \mu \| \partial_y \omega_k \|^2_{L^2} \le C_0 \delta_0 |k|^{-1} \operatorname{Dis}_{\omega, 4} + \delta^* \operatorname{Dis}_{\omega, 3} +  \frac{1}{4\delta^*} \mu^{-\frac{1}{3}} |k|^{\frac{4}{3}} \|  \theta_k\|^2_{L^2}.
	\end{align}
	\begin{proof}
		Taking the inner product of the first equation of system \eqref{linear_eq} with $ \omega_k$ in $L^2$, and then taking the real part $ \operatorname{Re}$, we obtain
		\begin{align}\label{1w_one}
		\frac{1}{2} \frac{d}{dt} \| \omega_k\|^2_{L^2} + \mu \| \partial_y \omega_k \|^2_{L^2} = \operatorname{Re} \langle ik U'' \psi_k, \omega_k \rangle - \operatorname{Re} \langle ik \theta_k, \omega_k \rangle.
		\end{align}
		For the non-local terms involving $ ik U'' \psi_k$ in the above expression, since $ \omega_k = \Delta_{k} \psi_k$, by performing integration by parts, we can directly obtain
		\begin{align}\label{ikU''_1}
			\big| \operatorname{Re} \langle ik U'' \psi_k, \omega_k \rangle \big| &= \big| \operatorname{Re} \langle ik U'' \nabla_k \psi_k, \nabla_k \psi_k \rangle + \operatorname{Re} \langle ik U''' ik \psi_k, \partial_y \psi_k \rangle \big|  \nn\\
			&\le |k|^2 \| U'''\|_{L^\infty} \| \psi_k\|_{L^2} \| \partial_{y} \psi_k\|_{L^2} \le C_0 \delta_0 |k|^{-1} \operatorname{Dis}_{\omega, 4},
		\end{align}
		where we have used the fact $ \operatorname{Re} \langle ik U'' \nabla_k \psi_k, \nabla_k \psi_k \rangle = 0.$ For the terms involving the non-local term $ ik \theta_k$, we directly apply Young's inequality to obtain
		\begin{align}\label{ikta_1}
			\big| \operatorname{Re} \langle ik \theta_k , \omega_k \rangle \big| \le \delta^* \operatorname{Dis}_{\omega, 3} +  \frac{1}{4\delta^*} \mu^{-\frac{1}{3}} |k|^{\frac{4}{3}} \|  \theta_k\|^2_{L^2}.
		\end{align}
		Collecting \eqref{1w_one}--\eqref{ikta_1} gives rise to \eqref{1w}.
	\end{proof}
\end{lemma}

\begin{lemma}\label{le_w2}
	Under the conditions of Theorem \ref{theorem1.1}, we have
	\begin{align}\label{2w}
	&\quad \frac{1}{2} \frac{d}{dt} \mu^{\frac{2}{3}} |k|^{-\frac{2}{3}}\| \partial_y \omega_k\|^2_{L^2} + \frac{1}{2}\mu^{\frac{5}{3}} |k|^{-\frac{2}{3}} \| \partial_{yy} \omega_k \|^2_{L^2} \nn\\
	& \le (16 + \frac{C_0}{2}\delta_0)\operatorname{Dis}_{\omega, 1} + (\frac{1}{16} + \frac{C_0}{2}\delta_0) \operatorname{Dis}_{\omega, 3} + C_0 \mu^{-\frac{1}{3}} |k|^{\frac{4}{3}} \|  \theta_k\|^2_{L^2}.
	\end{align}
	\begin{proof}
		By applying a similar procedure as in \eqref{1w_one}, we obtain
		\begin{align}\label{2w_one}
			\frac{1}{2} \frac{d}{dt} \mu^{\frac{2}{3}} |k|^{-\frac{2}{3}}\| \partial_y \omega_k\|^2_{L^2} + \mu^{\frac{5}{3}} |k|^{-\frac{2}{3}} \| \partial_{yy} \omega_k \|^2_{L^2} =- \mu^{\frac{2}{3}} |k|^{-\frac{2}{3}}\operatorname{Re} \langle ik\big( U'  \omega_k -  \partial_y (U'' \psi_k) +  \partial_y \theta_k \big), \partial_y \omega_k \rangle.
		\end{align}
		Here, we have used integration by parts and the fact that $$ \partial_{yy} \omega_k|_{y = \pm 1} = \int_{-1}^{1} \partial^2_{yy} G_{k}(y, y') \omega_k(y') \mathrm{d} y'|_{y = \pm 1} = 0.$$
		Next, we need to estimate each term on the right-hand side of the above expression. For the first term, using an estimate similar to that in \eqref{1ta_one}, we obtain
		\begin{align}\label{2w1}
			\mu^{\frac{2}{3}} |k|^{-\frac{2}{3}}\operatorname{Re} \langle ik U'  \omega_k  , \partial_y \omega_k \rangle \le  16\mu \| \partial_y \omega_k\|^2_{L^2} + \frac{1}{16}\mu^{\frac{1}{3}} |k|^{\frac{2}{3}} \|  \omega_k\|^2_{L^2}.
		\end{align}
		For the second term, we directly apply \eqref{W_estimate} to give rise to
		\begin{align}\label{2w2}
			\mu^{\frac{2}{3}} |k|^{-\frac{2}{3}}\operatorname{Re} \langle ik \partial_{y} (U''  \psi_k)  , \partial_y \omega_k \rangle &\lesssim \mu^{\frac{2}{3}} |k|^{-\frac{2}{3}} \big( \| U'''\|_{L^\infty} |k| \| \psi_k\|_{L^2} + \| U''\|_{L^\infty} |k| \| \partial_y \psi_k \|_{L^2}\big) \| \partial_y \omega_k\|_{L^2}\nn\\
			&\le C_0 \delta_0 \mu^{\frac{2}{3}} |k|^{-\frac{2}{3}} \| \omega_k\|_{L^2} \| \partial_y \omega_k\|_{L^2}\nn\\
			& \le \frac{C_0}{2} \delta_0  \mu \| \partial_y \omega_k\|^2_{L^2} + \frac{C_0}{2} \delta_0 \mu^{\frac{1}{3}} |k|^{\frac{2}{3}} \|  \omega_k\|^2_{L^2},
		\end{align}
		where we have used the following fact
		\begin{align}\label{w_fact}
			\| \omega_k\|^2_{L^2} = \int_{-1}^{1} (-k^2 \psi_k + \partial_{yy} \psi_k )^2 \mathrm{d} y &= k^4 \| \psi_k\|^2_{L^2} + \| \partial_{yy} \psi_k \|^2_{L^2} + 2|k|^2 \| \partial_{y} \psi_k\|^2_{L^2}\nn\\
			 &\ge \frac{1}{2} ( |k|\| \psi_k\|_{L^2} + |k|\| \partial_{y} \psi_k\|_{L^2} )^2.
		\end{align}
		For the third term, by performing integration by parts and then applying Young's inequality, we can get
		\begin{align}\label{2w3}
			\mu^{\frac{2}{3}} |k|^{-\frac{2}{3}}\operatorname{Re} \langle ik \partial_{y} \theta_k  , \partial_y \omega_k \rangle \le \frac{1}{2} \mu^{\frac{5}{3}} |k|^{-\frac{2}{3}} \| \partial_{yy} \omega_k \|^2_{L^2} + \frac{1}{2} \mu^{-\frac{1}{3}} |k|^{\frac{4}{3}} \|  \theta_k \|^2_{L^2}.
		\end{align}
		By substituting \eqref{2w1}, \eqref{2w2}, and \eqref{2w3} into \eqref{2w_one}, we obtain \eqref{2w}.
	\end{proof}
\end{lemma}
Next, we consider the following cross terms to induce the enhanced dissipation effect of $ \omega_k$.
\begin{lemma}\label{le_w3}
	Under the conditions of Theorem \ref{theorem1.1}, we have
	\begin{align}\label{3w}
	&\quad  \mu^{\frac{1}{3}} |k|^{-\frac{4}{3}}\frac{d}{dt} \operatorname{Re} \langle ik \omega_k, \partial_y \omega_k \rangle  + \mu^{\frac{1}{3}} |k|^{\frac{2}{3}} \| \sqrt{U'} \omega_k \|^2_{L^2}  \nn\\
	& \le  \frac{1}{2}\operatorname{Dis}_{\omega, 1} +  C_0\delta_0 \operatorname{Dis}_{\omega, 3} + \operatorname{Dis}_{\omega, 1} + \operatorname{Dis}_{\omega, 2}+ C_0 \mu^{-\frac{1}{3}} |k|^{\frac{4}{3}} \|  \theta_k\|^2_{L^2}.
	\end{align}
	\begin{proof}
		First, we compute the following time derivative, which gives
		\begin{align}\label{3w_one}
			\mu^{\frac{1}{3}} |k|^{-\frac{4}{3}} \frac{d}{dt} \operatorname{Re} \langle ik \omega_k, \partial_{y} \omega_k \rangle & = \mu^{\frac{1}{3}} |k|^{-\frac{4}{3}} \operatorname{Re} \langle ik(-ik U\omega_k + ik U'' \psi_k + \mu \partial_{yy} \omega_k - ik \theta_k), \partial_y \omega_k \rangle\nn\\
			& \quad + \mu^{\frac{1}{3}} |k|^{-\frac{4}{3}} \operatorname{Re} \langle ik \omega_k, \partial_y (-ik U\omega_k + ik U'' \psi_k + \mu \partial_{yy} \omega_k - ik \theta_k)  \rangle\nn\\
			&  = 2\mu^{\frac{1}{3}} |k|^{-\frac{4}{3}} \Big( \operatorname{Re} \langle |k|^2 U\omega_k, \partial_y \omega_k \rangle - \operatorname{Re} \langle |k|^2 U'' \psi_k, \partial_y \omega_k \rangle \nn\\
			& \quad + \mu\operatorname{Re} \langle ik\partial_{yy} \omega_k, \partial_y \omega_k \rangle + |k|^2 \operatorname{Re} \langle  \theta_k, \partial_y \omega_k \rangle \Big).
		\end{align}
		For the first term above, we can directly obtain the result by performing integration by parts, yielding
		\begin{align}\label{3w1}
			2\mu^{\frac{1}{3}} |k|^{-\frac{4}{3}} \operatorname{Re} \langle |k|^2 U\omega_k, \partial_y \omega_k \rangle  = - \mu^{\frac{1}{3}} |k|^{\frac{2}{3}} \| \sqrt{U'} \omega_k\|^2_{L^2}.
		\end{align}
		For the second term, by applying integration by parts and utilizing the fact from \eqref{w_fact}, we obtain
		\begin{align}\label{3w2}
			-2\mu^{\frac{1}{3}} |k|^{-\frac{4}{3}} \operatorname{Re} \langle |k|^2 U'' \psi_k, \partial_y \omega_k \rangle  &\le 2 \mu^{\frac{1}{3}} |k|^{-\frac{1}{3}} \big( \|U'' \|_{L^\infty} \| k \partial_{y}\psi_k\|_{L^2} + \|U''' \|_{L^\infty} \| k \psi_k\|_{L^2}\big)\| \omega_k\|_{L^2}\nn\\
			& \le C_0 \delta_0  \operatorname{Dis}_{\omega, 3}.
		\end{align}
		For the last two terms, we directly apply Young's inequality to obtain
		\begin{align}\label{3w3}
			&\quad 2\mu^{\frac{1}{3}} |k|^{-\frac{4}{3}} \Big( \mu\operatorname{Re} \langle ik\partial_{yy} \omega_k, \partial_y \omega_k \rangle +|k|^2 \operatorname{Re} \langle  \theta_k, \partial_y \omega_k \rangle\Big) \nn\\
			& \le \frac{1}{2} \operatorname{Dis}_{\omega, 1} + 2 \operatorname{Dis}_{\omega, 1}^{\frac{1}{2}} \operatorname{Dis}_{\omega, 2}^{\frac{1}{2}} + 2\mu^{-\frac{1}{3}} |k|^{\frac{4}{3}} \| \theta_k\|^2_{L^2}.
		\end{align}
		Collecting \eqref{3w_one}--\eqref{3w3} gives rise to \eqref{3w}.
	\end{proof}
\end{lemma}

\begin{lemma}\label{le_w4}
		Under the conditions of Theorem \ref{theorem1.1}, we have
	\begin{align}\label{4w}
	   \frac{d}{dt} \operatorname{Re} \langle  \omega_k, \mathfrak{J}_k [\omega_k] \rangle  + \frac{1}{4} \operatorname{Dis}_{\omega, 4}  \le  \frac{C_0}{4}\operatorname{Dis}_{\omega, 1} +  C_0 \operatorname{Dis}_{\omega, 2} + \delta^* \operatorname{Dis}_{\omega, 3}+ \frac{C_0^2}{4\delta^*} \mu^{-\frac{1}{3}} |k|^{\frac{4}{3}} \|  \theta_k\|^2_{L^2}.
	\end{align}
	\begin{proof}
		By directly computing the time derivative and applying Lemma \ref{le2.3}, we can get
		\begin{align}\label{4w_one}
			\frac{d}{dt} \operatorname{Re} \langle \omega_k, \mathfrak{J}_k[\omega_k] \rangle &=  \operatorname{Re} \int_{-1}^{1} \frac{d}{dt} \omega_k \cdot \overline{ \mathfrak{J}_k[\omega_k]} \mathrm{d} y + \operatorname{Re} \int_{-1}^{1}  \omega_k \cdot \overline{ \mathfrak{J}_k[\frac{d}{dt} \omega_k]} \mathrm{d} y\nn\\
			& = 2\operatorname{Re} \int_{-1}^{1} \frac{d}{dt} \omega_k \cdot \overline{ \mathfrak{J}_k[\omega_k]} \mathrm{d} y\nn\\
			& = 2\operatorname{Re} \langle \mu \partial_{yy} \omega_k + ik U'' \psi_k + ikU \omega_k   - ik \theta_k, \mathfrak{J}_{k} [\omega_k] \rangle\nn\\
			& := \sum_{i=1}^{4}H_i.
		\end{align}
		For the first term $H_1$, by applying Lemma \ref{le2.1} and Lemma \ref{le2.3}, we obtain
		\begin{align}\label{H1}
			H_1 & = 2\mu \operatorname{Re} \int_{-1}^{1} \partial_{yy} \omega_k \cdot \overline{\mathfrak{J}_k [\omega_k]}  \mathrm{d} y = 2 \mu \operatorname{Re} \int_{-1}^{1} \mathfrak{J}_k [\partial_{yy} \omega_k]   \cdot \overline{ \omega_k} \mathrm{d} y \nn\\
			& \lesssim \mu \| \partial_{yy} \omega_k \|_{L^2} \| \omega_k\|_{L^2}  \le \frac{C_0}{4} \operatorname{Dis}_{\omega, 1} +  C_0 \operatorname{Dis}_{\omega, 2}.
		\end{align}
		For the second term $ H_2$, by applying Lemma \ref{le2.1}-Lemma \ref{le2.2}, and the fact that $ \omega_k = (-k^2 + \partial_{yy}) \psi_k$, we can get
		\begin{align}\label{H2}
			H_2 & = 2 \operatorname{Re} \int_{-1}^{1} ik U'' \psi_k \cdot \overline{\mathfrak{J}_{k} [-k^2 \psi_k]} + 2 \operatorname{Re} \int_{-1}^{1} ik U'' \psi_k \cdot \overline{\mathfrak{J}_{k} [\partial_{yy} \psi_k]} \mathrm{d} y\nn\\
			& = - 2 \operatorname{Re} \int_{-1}^{1} (ik) ik U'' \psi_k \cdot \overline{\mathfrak{J}_{k} [ik \psi_k]} + \partial_{y} (ik U'' \psi_k) \cdot \overline{\mathfrak{J}_{k} [\partial_{y} \psi_k]} -  ik U'' \psi_k \cdot \overline{[\mathfrak{J}_{k}, \partial_{y}] [\partial_{y} \psi_k]} \mathrm{d} y\nn\\
			& \lesssim |k| \| (U'', U''')\|_{L^\infty} \big( \| \mathfrak{J}_k\|_{L^2} + |k|^{-1} \| [\mathfrak{J}_k, \partial_{y}]\|_{L^2} \big) \| \nabla_k \psi_k\|^2_{L^2}\nn\\
			& \le C_0 \delta_0 \| \nabla_k \psi_k\|^2_{L^2}.
		\end{align}
		For the fourth term $ H_4$, we can  apply Lemma \ref{le2.1} and Young's inequality to obtain
		\begin{align}\label{H4}
			H_4 \le C_0|k| \| \theta_k\|_{L^2} \| \omega_k\|_{L^2} \le \delta^* \operatorname{Dis}_{\omega, 3}+ \frac{C_0^2}{4\delta^*} \mu^{-\frac{1}{3}} |k|^{\frac{4}{3}} \|  \theta_k\|^2_{L^2}.
		\end{align}
		The most challenging term here is $H_3$. First, by applying Lemma \ref{le2.3}, we obtain
		\begin{align}\label{H3_split}
			H_3 &= \operatorname{Re} \Big( -ik \int_{-1}^{1} U \omega_k \cdot \overline{\mathfrak{J}_k [\omega_k]} \mathrm{d} y + k \int_{-1}^{1} iU \omega_k \cdot \mathfrak{J}_k [\overline{\omega_k}] \mathrm{d} y\Big)\nn\\
			& = \operatorname{Re} \Big( -ik \int_{-1}^{1} U \omega_k \cdot \overline{\mathfrak{J}_k [\omega_k]} \mathrm{d} y +   k \int_{-1}^{1} \mathfrak{J}_k [iU \omega_k] \cdot  \overline{\omega_k} \mathrm{d} y\Big)\nn\\
			& = \frac{k^2}{2} \operatorname{Re} \int_{-1}^{1} \int_{-1}^{1}  \frac{U(y)-U(y')}{y-y'} G_k(y,y') \overline{\omega_k(y)} \omega_k(y')   \mathrm{d} y' \mathrm{d} y.
		\end{align}
		Recalling the integral remainder form of Taylor's formula, we have
		 \begin{align*}
		 	\frac{U(y)-U(y')}{y-y'} = U'(y) + \frac{1}{2}U''(y) (y-y') + \frac{1}{6} h(y,y'),
		 \end{align*}
		 where  $ h(y,y') = \frac{1}{y-y'} \int_{y}^{y'} U'''(s) (y' - s)^2 \mathrm{d} s$. Then, term $H_3$ can be decomposed into the following three terms for further analysis
		 \begin{align*}
		 	H_3 = H_{31} + H_{32}  + H_{33}.
		 \end{align*}
		 Here, the right-hand side of the above equation is defined as follows:
		 \begin{align*}
		 	&H_{31} := \frac{k^2}{2} \operatorname{Re} \int_{-1}^{1} \int_{-1}^{1}  U'(y) G_k(y,y') \overline{\omega_k(y)} \omega_k(y')   \mathrm{d} y' \mathrm{d} y,\\
		 	&H_{32} := \frac{k^2}{4} \operatorname{Re} \int_{-1}^{1} \overline{\omega_k(y)} \int_{-1}^{1}  \mathscr{K}_1(y,y')   \omega_k(y')   \mathrm{d} y' \mathrm{d} y,\\
		 	&H_{33} := \frac{k^2}{12} \operatorname{Re} \int_{-1}^{1} \overline{\omega_k(y)} \int_{-1}^{1}  \mathscr{K}_2(y,y')  \omega_k(y')   \mathrm{d} y' \mathrm{d} y,
		 \end{align*}
		 where the kernel functions $ \mathscr{K}_1(y,y')$ and $ \mathscr{K}_2(y,y')$ are defined in the Appendix.
		 Since $G_k(y,y') $ and $ \psi_k(y')$ vanish at the boundary $ y' = \pm 1$, we can obtain the result by applying integration by parts
		 \begin{align}\label{H31}
		 	H_{31} &= \frac{k^2}{2} \operatorname{Re} \int_{-1}^{1}  \overline{\omega_k(y)} U'(y) \int_{-1}^{1} G_k(y,y')  \Delta_{k} \psi_k(y')   \mathrm{d} y' \mathrm{d} y \nn\\
		 	& = \frac{k^2}{2} \operatorname{Re} \int_{-1}^{1}    U'(y)  \psi_k(y)  \overline{\Delta_{k} \psi_k(y)}  \mathrm{d} y\nn\\
		 	& = -\frac{k^2}{2} \| \sqrt{U'} \nabla_k \psi_k\|^2_{L^2} - \frac{k^2}{2} \operatorname{Re} \int_{-1}^{1} U''(y)  \psi_k  \overline{ \partial_{y} \psi_k}  \mathrm{d} y\nn\\
		 	& \le -\frac{k^2}{2} \| \sqrt{U'} \nabla_k \psi_k\|^2_{L^2} + C_0 \delta_0 \operatorname{Dis}_{\omega, 4}.
		 \end{align}
		 For term $H_{32}$, applying integration by parts twice and using the boundary condition of $G_k(y,y')$ along with Lemma \ref{A2}, we obtain
		 \begin{align}\label{H32}
		 	H_{32} &= \frac{k^2}{4} \operatorname{Re} \int_{-1}^{1}  \overline{\nabla_k \psi_k(y)} \int_{-1}^{1} \nabla_k \big( ik \mathscr{K}_1(y,y') \cdot  ik\psi_k(y') + \partial_{y'} \mathscr{K}_1(y,y')  \partial_{y'} \psi_k(y')\big) \mathrm{d} y' \mathrm{d} y \nn\\
		 	& \lesssim k^2 \| \nabla_k \psi_k\|^2_{L^2}\big( k^2 \|\mathscr{K}_1(y,y')  \|_{L^2_{y,y'}} + |k| \|(\partial_{y}\mathscr{K}_1(y,y'), \partial_{y'}\mathscr{K}_1(y,y') ) \|_{L^2_{y,y'}} + \|\partial_{yy'}\mathscr{K}_1(y,y') \|_{L^2_{y,y'}} \big)\nn\\
		 	&\le C_0 \delta_0 \operatorname{Dis}_{\omega, 4}.
		 \end{align}
		 Similarly, it follows from Lemma \ref{A3} that
		 \begin{align}\label{H33}
		 	H_{33}  &= \frac{k^2}{12} \operatorname{Re} \int_{-1}^{1}  \overline{\nabla_k \psi_k(y)} \int_{-1}^{1} \nabla_k \big( ik \mathscr{K}_2(y,y') \cdot  ik\psi_k(y') + \partial_{y'} \mathscr{K}_2(y,y')  \partial_{y'} \psi_k(y')\big) \mathrm{d} y' \mathrm{d} y \nn\\
		 	& \lesssim k^2 \| \nabla_k \psi_k\|^2_{L^2}\big( k^2 \|\mathscr{K}_2(y,y')  \|_{L^2_{y,y'}} + |k| \|(\partial_{y}\mathscr{K}_2(y,y'), \partial_{y'}\mathscr{K}_2(y,y') ) \|_{L^2_{y,y'}} + \|\partial_{yy'}\mathscr{K}_2(y,y') \|_{L^2_{y,y'}} \big)\nn\\
		 	&\le C_0 \delta_0 \operatorname{Dis}_{\omega, 4}.
		 \end{align}
		 Collecting \eqref{H31}--\eqref{H33} gives rise to
		 \begin{align}\label{H3}
		 	H_3 \le -\frac{k^2}{2} \| \sqrt{U'} \nabla_k \psi_k\|^2_{L^2} + 3C_0 \delta_0 \operatorname{Dis}_{\omega, 4}.
		 \end{align}
		 Combining the estimates for $H_1$-$H_4$, we can immediately obtain \eqref{4w} when $\delta_0$ is sufficiently small.
	\end{proof}
\end{lemma}

\begin{lemma}\label{le_w5}
	Under the conditions of Theorem \ref{theorem1.1}, we have
	\begin{align}\label{5w}
		\frac{d}{dt}E_{\omega, k} + 2\operatorname{Dis}_{\omega,1} + 2 \operatorname{Dis}_{\omega,2} + \frac{1}{4}\operatorname{Dis}_{\omega,3} + \frac{1}{16}c_{\alpha}\operatorname{Dis}_{\omega,4} \le C_0 \mu^{-\frac{1}{3}} |k|^{\frac{4}{3}} \| \theta_k\|^2_{L^2}.
	\end{align}
	\begin{proof}
		First, setting $\delta^* = \frac{1}{16^3}$, then multiplying equation \eqref{1w} by 256, \eqref{2w} by 8, adding \eqref{3w}, and adding $c_{\alpha}$ times \eqref{4w}, we obtain
		\begin{align}
			&\quad \frac{d}{dt} E_{\omega, k} + 4\operatorname{Dis}_{\omega,1} + 2 \operatorname{Dis}_{\omega,2} + \frac{1}{2}\operatorname{Dis}_{\omega,3} + \frac{1}{8}c_{\alpha}\operatorname{Dis}_{\omega,4} \nn\\
			&\le 256 C_0\delta_0\operatorname{Dis}_{\omega, 4} + C_0c_{\alpha}\operatorname{Dis}_{\omega, 1} + \frac{c_{\alpha}}{16^3} \operatorname{Dis}_{\omega, 3} + C_0 \mu^{-\frac{1}{3}} |k|^{\frac{4}{3}} \| \theta_k\|^2_{L^2}.
		\end{align}
		Further setting $c_{\alpha} = \min\{ \frac{4}{C_0}, 1\}$, we obtain \eqref{5w}, which completes the proof of this lemma.
	\end{proof}
\end{lemma}
\subsection{Completion the proof of  Theorem \ref{theorem1.1}}
From the definition of $ E_{\theta, k}$, it follows that $ E_{\theta, k} \approx \| \theta_k\|^2_{L^2} + \nu^{\frac{2}{3}} |k|^{-\frac{2}{3}}\| \partial_{y} \theta_k\|^2_{L^2}$ is a coercive energy functional. Then, by applying Lemma \ref{le_ta1}, we obtain
\begin{align*}
\frac{d}{dt} E_{\theta, k} + \frac{1}{4} \nu^{\frac{1}{3}} |k|^{\frac{2}{3}} \big(\|  \theta_k\|^2_{L^2} + \nu^{\frac{2}{3}} |k|^{-\frac{2}{3}}\| \partial_{y} \theta_k\|^2_{L^2}\big)  \le 0.
\end{align*}
Applying Gronwall's inequality to the above equation, we directly obtain the result of \eqref{linear_ta} in Theorem \ref{theorem1.1}. Similarly, $ E_{\omega, k} \approx \| \omega_k\|^2_{L^2} + \mu^{\frac{2}{3}} |k|^{-\frac{2}{3}}\| \partial_{y} \omega_k\|^2_{L^2}$ is also a coercive energy functional. We can simplify \eqref{5w} as follows,
\begin{align}
	\frac{d}{dt}E_{\omega, k} + \frac{1}{4} \mu^{\frac{1}{3}} |k|^{\frac{2}{3}} \big(\|  \| \omega_k\|^2_{L^2} + \mu^{\frac{2}{3}} |k|^{-\frac{2}{3}}\| \partial_{y} \omega_k\|^2_{L^2}\big)  \le C_0 \mu^{-\frac{1}{3}} |k|^{\frac{4}{3}} \| \theta_k\|^2_{L^2}.
\end{align} and then applying Gronwall's inequality, we also obtain the result of \eqref{linear_w}. This completes the proof of Theorem \ref{theorem1.1}.
\hspace{12.3cm}$\square$

\section{Nonlinear stability}
In this section, we focus on the case $\sigma = 1$ to establish the proof of Theorem \ref{thm1.2}. The proof is divided into four subsections. The first subsection introduces some energy functionals and preliminary lemmas. The second subsection presents the nonlinear estimates for the zero mode. The third subsection provides the nonlinear estimates for the nonzero modes. The fourth subsection employs the bootstrap method to establish global nonlinear stability.

\subsection{Notations and Preliminary Lemmas}For convenience, we slightly modify the dissipative energy functional from the linear stability analysis in  Section \ref{sec2} as follows:
\begin{align*}
&\operatorname{Dis}_{\theta, 1} = \nu \| \nabla_k \theta_k\|^2_{L^2}, \quad \operatorname{Dis}_{\theta, 2} = \nu^{\frac{5}{3}} |k|^{-\frac{2}{3}} \| \partial_{y} \nabla_k \theta_k\|^2_{L^2},\quad
\operatorname{Dis}_{\theta, 3} = \nu^{\frac{1}{3}} |k|^{\frac{2}{3}} \|  \theta_k\|^2_{L^2},
\end{align*}
and
\begin{align*}
&\operatorname{Dis}_{\omega, 1} = \mu \| \nabla_k \omega_k\|^2_{L^2}, \quad \operatorname{Dis}_{\omega, 2} = \mu^{\frac{5}{3}} |k|^{-\frac{2}{3}} \| \partial_{y} \nabla_k \omega_k\|^2_{L^2},\quad\operatorname{Dis}_{\omega, 3} = \mu^{\frac{1}{3}} |k|^{\frac{2}{3}} \|  \omega_k\|^2_{L^2},\nn\\
& \operatorname{Dis}_{\omega, 4} = |k|^2 \| \nabla_k \psi_k\|^2_{L^2},
\quad\operatorname{Dis}_{\omega, 5} = \mu^{\frac{2}{3}} |k|^{ \frac{4}{3}} \| \partial_{y}\nabla_k \psi_k\|^2_{L^2} \quad (\text{where} \nabla_k = (ik, \partial_{y})^\top ).
\end{align*}
Additionally, for the nonlinear stability analysis, we first set $ \lambda = \min\{\mu, \nu\}$, and then we also need to introduce the energy functional corresponding to the linear stability as follows:
\begin{align*}
\mathscr{E}_{\theta} =& 	\mathscr{E}_{\theta,0} + \mathscr{E}_{\theta,\neq},
\quad\hbox{\it with}\quad\\
&\mathscr{E}_{\theta,0} = 16 e^{2\delta_1 \lambda t} \| \theta_0 \|^2_{L^2} +   e^{2\delta_1 \lambda t} \nu^{\frac{2}{3}} \| \partial_{y} \theta_0 \|^2_{L^2}, \quad
\mathscr{E}_{\theta,\neq} = e^{2\delta_1 \lambda^{\frac{1}{3}} t} \sum_{k\neq 0} |k|^{2m} E_{\theta, k}[\theta_k].
\\
\mathscr{E}_{\omega} =& 	\mathscr{E}_{\omega,0} + \mathscr{E}_{\omega,\neq}
,\quad\hbox{\it with}\quad\\
&\mathscr{E}_{\omega,0} = 128 e^{2\delta_1 \lambda t} \| \omega_0 \|^2_{L^2} + 4 e^{2\delta_1 \lambda t} \mu^{\frac{2}{3}} \| \partial_{y} \omega_0 \|^2_{L^2}, \quad
\mathscr{E}_{\omega,\neq} = e^{2\delta_1 \lambda^{\frac{1}{3}} t}\sum_{k\neq 0} |k|^{2m} E_{\omega, k}[\omega_k].
\\
 \mathscr{D}_{\theta} =& \mathscr{D}_{\theta, 0} + \mathscr{D}_{\theta, \neq}= \mathscr{D}_{\theta, 0} +\sum_{i=1}^{3} \mathscr{D}_{\theta, i},\quad\hbox{\it with}\quad\\
&\mathscr{D}_{\theta, 0} =  32 e^{2\delta_1 \lambda t} \nu \| \partial_{y} \theta_0 \|^2_{L^2} +   2 e^{2\delta_1 \lambda t} \nu^{\frac{5}{3}}  \| \partial_{yy} \theta_0 \|^2_{L^2}, \quad \mathscr{D}_{\theta, i} = e^{2\delta_1 \lambda^{\frac{1}{3}} t} \sum_{k\neq 0} |k|^{2m} \operatorname{Dis}_{\theta, i}.
\\
\mathscr{D}_{\omega} =& \mathscr{D}_{\omega, 0} + \mathscr{D}_{\omega, \neq}
= \mathscr{D}_{\omega, 0} + \sum_{i=1}^{5} \mathscr{D}_{\omega, i},\quad\hbox{\it with}\quad\\
& \mathscr{D}_{\omega, 0} = 256 e^{2\delta_1 \lambda t}  \mu \| \partial_{y} \omega_0 \|^2_{L^2} + 8 e^{2\delta_1 \lambda t}  \mu^{\frac{5}{3}}  \| \partial_{yy} \omega_0 \|^2_{L^2}, \quad \mathscr{D}_{\omega, i} = e^{2\delta_1 \lambda^{\frac{1}{3}} t} \sum_{k\neq 0} |k|^{2m} \operatorname{Dis}_{\omega, i}.
\end{align*}
Here, the coercive functional $ E_{\omega, k}[\omega_k]$ defined in Section \ref{sec2} is also modified as follows
\begin{align*}
	E_{\omega, k}[\omega_k] &= 128\| \omega_k\|^2_{L^2} + 4 \mu^{\frac{2}{3}} |k|^{-\frac{2}{3}} \|\partial_y \omega_k \|^2_{L^2} +  \mu^{\frac{1}{3}} |k|^{-\frac{4}{3}}
	\operatorname{Re} \langle i k \omega_k, \partial_y \omega_k \rangle\nn\\
	&\quad  + c_{\alpha} \operatorname{Re} \langle  \omega_k, \mathfrak{J}_k [\omega_k] \rangle +   c_{\beta} \mu^{\frac{2}{3}}|k|^{-\frac{2}{3}} \operatorname{Re} \langle  \partial_{y}\omega_k, \mathfrak{J}_k [\partial_{y} \omega_k] \rangle\nn\\
	& := E^1_{\omega, k}[\omega_k] +  c_{\beta} \mu^{\frac{2}{3}}|k|^{-\frac{2}{3}} \operatorname{Re} \langle  \partial_{y}\omega_k, \mathfrak{J}_k [\partial_{y} \omega_k] \rangle.
\end{align*}
Unlike the case $\sigma = 0$ in Theorem \ref{theorem1.1}, when $\sigma = 1$, we can further obtain the dissipation of higher-order derivatives of the stream function $ \psi_k$. Specifically, we first analyze the linear stability of System \eqref{linear_eq}, we obtain
\begin{lemma}\label{le_w6}
	Under the conditions of Theorem \ref{thm1.2}, consider the stability of System \eqref{linear_eq}. Then for any integer $k \neq 0$ and any $t>0$, we have
		\begin{align}\label{6w}
		\frac{d}{dt}E_{\omega, k} +  \operatorname{Dis}_{\omega,1} +   \operatorname{Dis}_{\omega,2} + \frac{1}{8}\operatorname{Dis}_{\omega,3} + \frac{1}{16}c_{\alpha}\operatorname{Dis}_{\omega,4} + \frac{1}{4}c_{\beta}\operatorname{Dis}_{\omega,5} \le C_0 \mu^{-\frac{1}{3}} |k|^{\frac{4}{3}} \| \theta_k\|^2_{L^2},
		\end{align}
		and
		\begin{align}\label{6ta}
		\frac{d}{dt} E_{\theta, k} + \operatorname{Dis}_{\theta,1} +  2\operatorname{Dis}_{\theta,2}  + \frac{1}{4} \operatorname{Dis}_{\theta,3} \le 0.
		\end{align}
		\begin{proof}
			We first consider the higher-order derivative terms of the stream function $\psi_k$ in $E_{\omega, k}[\omega_k]$, by directly computing the following derivative, we obtain
			\begin{align*}
				\partial_{t} \partial_{y} \omega_k &=  \big( \mu \Delta_{k} - ikU + ik U'' \Delta_{k}^{-1}\big) \partial_{y} \omega_k  -ikU'\omega_k + ikU''\Delta_{k}^{-1}\omega_k - ik\partial_{y} \theta_k + ikU''[\partial_{y}, \Delta_{k}^{-1}] \omega_k \nn\\
				& := \sum_{i=1}^{7}R_{k,i}.
			\end{align*}
			Taking the direct $L^2$	inner product in the above equation, we obtain
			\begin{align}\label{6w1}
				\frac{d}{dt} \mu^{\frac{2}{3}}|k|^{-\frac{2}{3}} \operatorname{Re} \langle  \partial_{y}\omega_k, \mathfrak{J}_k [\partial_{y} \omega_k] \rangle = \sum_{i=1}^{7}  \mu^{\frac{2}{3}}|k|^{-\frac{2}{3}} \operatorname{Re} \langle  \partial_{y}\omega_k, \mathfrak{J}_k [R_{k,i} ] \rangle.
			\end{align}
			The treatment of $ R_{k,2} $ and $ R_{k,3} $ is similar to that of $H_2$ and $H_3$ in Section \ref{sec2}, and we can directly obtain
			\begin{align}\label{6w2}
				\sum_{i=2}^{ 3}  \mu^{\frac{2}{3}}|k|^{-\frac{2}{3}} \operatorname{Re} \langle  \partial_{y}\omega_k, \mathfrak{J}_k [R_{k,i} ] \rangle  \le - \frac{1}{2} \mu^{\frac{2}{3}}|k|^{\frac{5}{3}}  \| \sqrt{U'} \partial_{y} \nabla_k \psi_k\|^2_{L^2} +C_0 \delta_0 \operatorname{Dis}_{\omega,5}.
			\end{align}
			For the treatment of $R_{k,1}$, we need to apply Lemmas \ref{le2.1} to \ref{le2.3}, obtaining
			\begin{align}\label{Rk1}
				\mu^{\frac{2}{3}}|k|^{-\frac{2}{3}} \operatorname{Re} \langle  \partial_{y}\omega_k, \mathfrak{J}_k [R_{k,1} ] \rangle  &= \mu^{\frac{5}{3}}|k|^{-\frac{2}{3}} \operatorname{Re} \int_{-1}^{1} \mathfrak{J}_k[\Delta_{k} \partial_{y} \omega_k] \overline{ \partial_{y} \omega_k}  \mathrm{d} y\nn\\
				& =  \mu^{\frac{5}{3}}|k|^{-\frac{2}{3}} \operatorname{Re} \int_{-1}^{1}  \Big( -\mathfrak{J}_k [\partial_{y} \nabla_k \omega_k] \overline{ \partial_{y} \nabla_k \omega_k} + [\mathfrak{J}_k, \partial_{y}] \overline{\partial_{y} \omega_k}\Big )  \mathrm{d} y\nn\\
				& \le C_0 \operatorname{Dis}_{\omega, 2}.
			\end{align}
			For the three terms from $ R_{k,4}$ to $ R_{k,6}$, it follows from Lemma \ref{le2.2} that
			\begin{align}\label{Rk456}
				\sum_{i=4}^{ 6}  \mu^{\frac{2}{3}}|k|^{-\frac{2}{3}} \operatorname{Re} \langle  \partial_{y}\omega_k, \mathfrak{J}_k [R_{k,i} ] \rangle &\le C_0  \mu^{\frac{2}{3}}|k|^{\frac{1}{3}}  \| \partial_{y} \omega_k\|_{L^2} \| (\omega_k, \theta_k)\|_{L^2}\nn\\
				& \le \frac{ C_0}{2} \operatorname{Dis}_{\omega, 1}   +C_0 \operatorname{Dis}_{\omega, 3} + C_0\mu^{-\frac{1}{3}}|k|^{\frac{4}{3}} \| \theta_k\|^2_{L^2}.
			\end{align}
			For the term $ R_{k,7}$, using the fact that $ \| (\partial_{y'} G_{k}(y,y'), \partial_{y} G_{k}(y,y'))\|_{L^\infty_{y'y}} \lesssim 1$, we apply Young's inequality to obtain
			\begin{align}\label{Rk7}
				\mu^{\frac{2}{3}}|k|^{-\frac{2}{3}} \operatorname{Re} \langle  \partial_{y}\omega_k, \mathfrak{J}_k [R_{k,7} ] \rangle &\le C_0 \delta_0 \mu^{\frac{2}{3}}|k|^{\frac{1}{3}} \| \partial_{y} \omega_k \|_{L^2} \| [\partial_{y} , \Delta_{k}^{-1}]\omega_k\|_{L^2}\nn\\
				& \le  C_0 \delta_0 \mu^{\frac{2}{3}}|k|^{\frac{1}{3}} \| \omega_k \|_{L^2} \| \partial_{y} \omega_k \|_{L^2}\nn\\
				&\le \frac{ C_0}{2} \operatorname{Dis}_{\omega, 1}   + C_0 \operatorname{Dis}_{\omega, 3}.
			\end{align}
			Collecting \eqref{6w1}--\eqref{Rk7} gives rise to
			\begin{align}\label{6w3}
			&\quad \frac{d}{dt} \mu^{\frac{2}{3}}|k|^{-\frac{2}{3}} \operatorname{Re} \langle  \partial_{y}\omega_k, \mathfrak{J}_k [\partial_{y} \omega_k] \rangle + \frac{1}{4} \operatorname{Dis}_{\omega,5}\nn\\
			& \le  C_0 \operatorname{Dis}_{\omega, 1}  + C_0 \operatorname{Dis}_{\omega, 2}   +2C_0 \operatorname{Dis}_{\omega, 3} + C_0\mu^{-\frac{1}{3}}|k|^{\frac{4}{3}} \| \theta_k\|^2_{L^2}.
			\end{align}
			For the energy functional $ E^1_{\omega, k}[\omega_k]$ in $ E_{\omega, k}[\omega_k]$, we can obtain the result using a method similar to that of \eqref{5w}
			\begin{align}\label{6w4}
				\frac{d}{dt}E^1_{\omega, k} + 2\operatorname{Dis}_{\omega,1} + 2 \operatorname{Dis}_{\omega,2} + \frac{1}{4}\operatorname{Dis}_{\omega,3} + \frac{1}{16}c_{\alpha}\operatorname{Dis}_{\omega,4} \le C_0 \mu^{-\frac{1}{3}} |k|^{\frac{4}{3}} \| \theta_k\|^2_{L^2}.
			\end{align}
			Taking $ c_{\beta} = \min\{ \frac{1}{16C_0}, 1\}$ and combining \eqref{6w3} and \eqref{6w4}, we obtain \eqref{6w}. Furthermore, by following the proof of \eqref{1ta}, we derive \eqref{6ta},  which completes the proof of the lemma.
		\end{proof}
\end{lemma}

\subsection{Nonlinear Estimates For Zero Mode} In this subsection, we analyze the nonlinear estimates for the zero mode of System \eqref{rewrite1}. Specifically, we consider the following system:
\begin{eqnarray}\label{nonlinear_eq_zero}
\left\{\begin{aligned}
&\partial_t \omega_0  = \mu \partial_{yy}\omega_0 - (\nabla^{\top} \psi \cdot \nabla \omega)_0 ,\\
& \partial_t \theta_0 =  \nu \partial_{yy}\theta_0  - (\nabla^{\top} \psi \cdot \nabla \theta)_0.
\end{aligned}\right.
\end{eqnarray}
\begin{lemma}\label{le_w7}
	Under the conditions of Theorem \ref{thm1.2}, we have
	\begin{align}\label{7wta1}
			\frac{d}{dt} \mathcal{E}_{\theta,0} +\frac{1}{2} \mathscr{D}_{\theta, 0} + (c_0\nu -2\delta_1 \lambda)\mathcal{E}_{\theta,0} \le C_0 \nu^{-\frac{1}{2}} \mathscr{D}^{\frac{1}{2}}_{\theta} \mathscr{D}^{\frac{1}{2}}_{\omega} \mathcal{E}^{\frac{1}{2}}_{\theta} + C_0 \nu^{-\frac{1}{3}} \mathscr{D}_{\theta} \mathcal{E}^{\frac{1}{2}}_{\omega},
	\end{align}
	and
	\begin{align}\label{7wta2}
		\frac{d}{dt} \mathcal{E}_{\omega,0} +\frac{1}{2} \mathscr{D}_{\omega, 0} + (c_0\mu  -2\delta_1 \lambda) \mathcal{E}_{\omega,0} \le C_0 \mu^{-\frac{1}{2}} \mathscr{D}_{\omega} \mathcal{E}^{\frac{1}{2}}_{\omega}.
	\end{align}
	\begin{proof}
		By applying a direct $L^2$	energy estimate, we obtain
		\begin{align}\label{7w1}
		\frac{d}{dt} \mathcal{E}_{\theta,0}  = \mathcal{L}_{\theta, 0} + \mathcal{N}_{\theta, 0} + 2\delta_1 \lambda \mathcal{E}_{\theta,0},
		\end{align}
		where
		\begin{align*}
		&\mathcal{L}_{\theta, 0} = 32 e^{2\delta_1 \lambda t} \nu \operatorname{Re} \langle \partial_{yy} \theta_0, \theta_0  \rangle + 2 e^{2\delta_1 \lambda t} \nu^{\frac{5}{3}} \operatorname{Re} \langle \partial_{yyy} \theta_0, \partial_{y}\theta_0  \rangle\nn\\
		&\mathcal{N}_{\theta, 0} = -32 e^{2\delta_1 \lambda t} \operatorname{Re} \langle \theta_0,  (\nabla^{\top} \psi \cdot \nabla \theta)_0 \rangle - 2 e^{2\delta_1 \lambda t} \nu^{\frac{2}{3}} \operatorname{Re} \langle \partial_{y}\theta_0, \partial_{y} (\nabla^{\top} \psi \cdot \nabla \theta)_0     \rangle.
		\end{align*}
		For the linear term $ \mathcal{L}_{\theta, 0}$, the boundary conditions of $ \theta_0|_{y = \pm 1} = 0$ imply that $\theta_0$ and $ \partial_y \theta_0$ satisfy the Poincar\'e inequality ( $c_0$ denotes the Poincar\'e constant). Thus, we have
		\begin{align}\label{7w2}
			\mathcal{L}_{\theta, 0} \le  - \frac{1}{2} \mathscr{D}_{\theta, 0} - c_0\nu \mathcal{E}_{\theta,0}.
		\end{align}
	For the nonlinear term $ 	\mathcal{N}_{\theta, 0}$, we first perform the Fourier expansion of the zero mode of the advection term $\nabla^\top \psi \cdot \nabla \theta$, obtaining
		\begin{align}\label{7w3}
			\mathcal{N}_{\theta, 0} &= -32 e^{2\delta_1 \lambda t} \operatorname{Re} \langle \theta_0, \sum_{\ell \in \Z} \partial_y \psi_\ell \cdot i(-\ell)  \theta_{-\ell} - \sum_{\ell \in \Z} i\ell \psi_\ell \cdot \partial_{y} \theta_{-\ell}  \rangle\nn\\
			&\quad  -  2 e^{2\delta_1 \lambda t} \nu^{\frac{2}{3}} \operatorname{Re} \langle \partial_{y}\theta_0, \partial_{y} \sum_{\ell \in \Z} \partial_y \psi_\ell \cdot (-\ell)  \theta_{-\ell} - \sum_{\ell \in \Z} i\ell \psi_\ell \cdot \partial_{y} \theta_{-\ell}      \rangle\nn\\
			& = 32 e^{2\delta_1 \lambda t} \operatorname{Re} \langle \theta_0, \partial_y \big( \sum_{\ell \in \Z\setminus\{0\}}  i\ell\psi_\ell    \theta_{-\ell}  \big) \rangle  +2 e^{2\delta_1 \lambda t} \nu^{\frac{2}{3}}\operatorname{Re} \langle \partial_{y} \theta_0, \partial^2_y \big( \sum_{\ell \in \Z\setminus\{0\}}  i\ell\psi_\ell    \theta_{-\ell}  \big) \rangle\nn\\
			& := N_0 + N_1.
		\end{align}
		For the term $ N_0 $, applying integration by parts and the Gagliardo–Nirenberg inequality yields
		\begin{align*}
			|N_0| &= | 32 e^{2\delta_1 \lambda t} \operatorname{Re} \langle \partial_y \theta_0,   \sum_{\ell \in \Z\setminus\{0\}}  i\ell\psi_\ell    \theta_{-\ell}   \rangle| \lesssim e^{2\delta_1 \lambda t} \| \partial_y \theta_0\|_{L^2} \sum_{\ell \in \Z\setminus\{0\}} \| i\ell\psi_\ell\|_{L^\infty} \| \theta_{-\ell}\|_{L^2}\nn\\
			& \lesssim  e^{\delta_1 \lambda t} \nu^{-\frac{1}{2}} \mathscr{D}^{\frac{1}{2}}_{\theta, 0} \sum_{\ell \in \Z\setminus\{0\}} |\ell| \| \psi_\ell\|^{\frac{1}{2}}_{L^2} \| \partial_{y} \psi_\ell\|^{\frac{1}{2}}_{L^2}  \| \theta_{-\ell}\|_{L^2}\nn\\
& \lesssim e^{2\delta_1 \lambda t} \nu^{-\frac{1}{2}} \mathscr{D}^{\frac{1}{2}}_{\theta, 0} \sum_{\ell \in \Z\setminus\{0\}} |\ell| \| \nabla_k \psi_\ell\|_{L^2}   \| \theta_{-\ell}\|_{L^2}\nn\\
			& \lesssim \nu^{-\frac{1}{2}} \mathscr{D}^{\frac{1}{2}}_{\theta, 0} \mathscr{D}^{\frac{1}{2}}_{\omega, 4} \mathcal{E}^{\frac{1}{2}}_{\theta, \neq}.
		\end{align*}
		For the term $ N_{1}$, we can decompose it into the following two parts
		\begin{align*}
				|N_1| &= \Big| 2 e^{2\delta_1 \lambda t} \nu^{\frac{2}{3}}\operatorname{Re} \Big\langle \partial_{yy} \theta_0, \partial_y \big( \sum_{\ell \in \Z\setminus\{0\}}  i\ell\psi_\ell    \theta_{-\ell}  \big) \Big\rangle\Big|\nn\\
				& \lesssim e^{2\delta_1 \lambda t} \nu^{\frac{2}{3}} \| \partial_{yy} \theta_0\|_{L^2} \Big( \sum_{\ell \in \Z\setminus\{0\}} \| i\ell\partial_{y}\psi_\ell\|_{L^\infty} \| \theta_{-\ell}\|_{L^2} + \| i\ell\psi_\ell\|_{L^\infty} \| \partial_{y}\theta_{-\ell}\|_{L^2}\Big)\nn\\
				& := N_{11} + N_{12}.
		\end{align*}
		For the first term above, we can directly obtain the result applying \eqref{w_fact} and the Gagliardo–Nirenberg inequality
		\begin{align*}
			N_{11} &\lesssim e^{2\delta_1 \lambda t} \nu^{\frac{2}{3}} \| \partial_{yy} \theta_0\|_{L^2} \sum_{\ell \in \Z\setminus\{0\}} |\ell| \| \partial_{y}\psi_\ell\|^{\frac{1}{2}}_{L^2} \| \partial_{yy}\psi_\ell\|^{\frac{1}{2}}_{L^2} \| \theta_{-\ell}\|_{L^2} \nn\\
			& \lesssim e^{\delta_1 \lambda t} \nu^{\frac{2}{3}} \| \partial_{yy} \theta_0\|_{L^2} \cdot e^{\delta_1 \lambda t}\sum_{\ell \in \Z\setminus\{0\}} |\ell| \| \omega_\ell\|_{L^2} \| \theta_{-\ell}\|_{L^2}\nn\\
			& \lesssim \nu^{-\frac{1}{3}} \mathscr{D}^{\frac{1}{2}}_{\theta, 0} \mathscr{D}^{\frac{1}{2}}_{\theta, 3} \mathcal{E}^{\frac{1}{2}}_{\omega, \neq}.
		\end{align*}
		Similarity, we have
		\begin{align}\label{N12}
			N_{12} &\lesssim e^{2\delta_1 \lambda t} \nu^{\frac{2}{3}} \| \partial_{yy} \theta_0\|_{L^2} \sum_{\ell \in \Z\setminus\{0\}} |\ell| \| \partial_{y}\psi_\ell\|^{\frac{1}{2}}_{L^2} \|  \psi_\ell\|^{\frac{1}{2}}_{L^2} \| \partial_{y} \theta_{-\ell}\|_{L^2} \nn\\
			& \lesssim e^{2\delta_1 \lambda t} \nu^{\frac{2}{3}} \| \partial_{yy} \theta_0\|_{L^2} \sum_{\ell \in \Z\setminus\{0\}} |\ell|^{\frac{1}{2}} \| \nabla_\ell \theta_\ell\|_{L^2} \| \partial_{y} \theta_{-\ell}\|_{L^2}\nn\\
			& \lesssim \nu^{-\frac{1}{3}} \mathscr{D}^{\frac{1}{2}}_{\theta, 0} \mathscr{D}^{\frac{1}{2}}_{\omega, 4} \mathcal{E}^{\frac{1}{2}}_{\theta, \neq}.
		\end{align}
		Collecting \eqref{7w2}-\eqref{N12}, we obtain \eqref{7wta1} in the lemma. Next, we consider the nonlinear estimates for the zero mode of $\omega$.
		\begin{align}\label{7w7}
			\frac{d}{dt} \mathcal{E}_{\omega,0}  = \mathcal{L}_{\omega, 0} + \mathcal{N}_{\omega, 0} + 2\delta_1 \lambda \mathcal{E}_{\omega,0} ,
		\end{align}
		where
		\begin{align*}
		&\mathcal{L}_{\omega, 0} = 256 e^{2\delta_1 \lambda t} \nu \operatorname{Re} \langle \partial_{yy} \omega_0, \omega_0  \rangle + 8 e^{2\delta_1 \lambda t} \nu^{\frac{5}{3}} \operatorname{Re} \langle \partial_{yyy} \omega_0, \partial_{y}\omega_0  \rangle,\nn\\
		&\mathcal{N}_{\omega, 0} = -256 e^{2\delta_1 \lambda t} \operatorname{Re} \langle \omega_0,  (\nabla^{\top} \psi \cdot \nabla \omega)_0 \rangle - 8 e^{2\delta_1 \lambda t} \nu^{\frac{2}{3}} \operatorname{Re} \langle \partial_{y}\omega_0, \partial_{y} (\nabla^{\top} \psi \cdot \nabla \omega)_0     \rangle.
		\end{align*}
		For the linear term $ \mathcal{L}_{\omega, 0}$, the boundary conditions of $ \omega_0|_{y = \pm 1} = 0$ imply that $\omega_0$ and $ \partial_y \omega_0$ satisfy the Poincar\'e inequality ( $c_0$ denotes the Poincar\'e constant). We can get
		\begin{align*}
		\mathcal{L}_{\omega, 0} \le  - \frac{1}{2} \mathscr{D}_{\omega, 0} - c_0\mu \mathcal{E}_{\omega,0}.
		\end{align*}
		In fact, for the nonlinear term $(\u \cdot \nabla \omega)_0$, a similar treatment of the advection term $ (\u \cdot \nabla \theta)_0$ allows us to obtain
		\begin{align}\label{Nw0}
			\mathcal{N}_{\omega, 0} \lesssim \mu^{-\frac{1}{2}} \mathscr{D}_{\omega} \mathcal{E}^{\frac{1}{2}}_{\omega}.
		\end{align}
		Combining \eqref{7w7}-\eqref{Nw0}, we obtain \eqref{7wta2}, which completes the proof of the lemma.
	\end{proof}
	
\end{lemma}

\subsection{Nonlinear Estimates For Nonzero Mode} In this subsection, we analyze the nonlinear estimates for the nonzero mode of System \eqref{rewrite1}. Specifically, we consider the following system:
\begin{eqnarray}\label{nonlinear_eq}
\left\{\begin{aligned}
&\partial_t \omega_k  = \mathcal{L}_{\omega, k} + \mathcal{N}_{\omega, k} ,\\
& \partial_t \theta_k =  \mathcal{L}_{\theta, k} + \mathcal{N}_{\theta, k},
\end{aligned}\right.
\end{eqnarray}
where
\begin{align*}
	&\mathcal{L}_{\omega, k} = -ikU \omega_k + ikU'' \psi_k + \mu\Delta_{k}\omega_k - ik \theta_k, \quad \mathcal{L}_{\theta, k} = -ik U\theta_k + \nu \Delta_{k}\theta_k,\nn\\
	& \mathcal{N}_{\omega, k} = -(\nabla^{\top} \psi \cdot \nabla \omega)_k, \quad  \mathcal{N}_{\theta, k} = -(\nabla^{\top} \psi \cdot \nabla \theta)_k.
\end{align*}
According to the definitions of $ \mathscr{E}_{\theta,\neq}$ and $ \mathscr{E}_{\omega,\neq}$, we directly compute the following time derivative
\begin{align*}
\frac{d}{dt} \mathscr{E}_{\theta,\neq} =& e^{2\delta_1 \lambda^{\frac{1}{3}} t} \sum_{k \in \Z\setminus\{0\}} |k|^{2m} \frac{d}{dt} E_{\theta, k} + 2\delta_1 \lambda^{\frac{1}{3}} \mathscr{E}_{\theta,\neq}  = \mathcal{L}_{\theta, \neq} + \mathcal{N}_{\theta, \neq} + 2\delta_1 \lambda^{\frac{1}{3}} \mathscr{E}_{\theta,\neq},
\\
\\
\frac{d}{dt} \mathscr{E}_{\omega,\neq} =& e^{2\delta_1 \lambda^{\frac{1}{3}} t} \sum_{k \in \Z\setminus\{0\}} |k|^{2m} \frac{d}{dt} E_{\omega, k} + 2\delta_1 \lambda^{\frac{1}{3}} \mathscr{E}_{\omega,\neq} = \mathcal{L}_{\omega, \neq} + \mathcal{N}_{\omega, \neq}  + 2\delta_1 \lambda^{\frac{1}{3}} \mathscr{E}_{\omega,\neq},
\end{align*}
where
\begin{align*}
\mathcal{L}_{\theta, \neq} &= 32 e^{2\delta_1 \lambda^{\frac{1}{3}} t} \sum_{k \in \Z\setminus\{0\}} |k|^{2m} \operatorname{Re} \langle \theta_k,  \mathcal{L}_{\theta, k}\rangle + 2 e^{2\delta_1 \lambda^{\frac{1}{3}} t} \nu^{\frac{2}{3}} \sum_{k \in \Z\setminus\{0\}} |k|^{2m-\frac{2}{3}} \operatorname{Re} \langle \partial_{y} \theta_k,  \partial_{y} \mathcal{L}_{\theta, k}\rangle\nn\\
&\quad  +  e^{2\delta_1 \lambda^{\frac{1}{3}} t} \nu^{\frac{1}{3}}  \sum_{k \in \Z\setminus\{0\}} |k|^{2m-\frac{4}{3}} \Big[ \operatorname{Re} \langle ik  \theta_k,  \partial_{y} \mathcal{L}_{\theta, k}\rangle  + \operatorname{Re} \langle ik \mathcal{L}_{\theta, k},  \partial_{y}\theta_k \rangle\Big],\nn\\
\mathcal{L}_{\omega, \neq} &= 256 e^{2\delta_1 \lambda^{\frac{1}{3}} t} \sum_{k \in \Z\setminus\{0\}} |k|^{2m} \operatorname{Re} \langle \omega_k,  ( 1+ \frac{c_{\alpha}}{128} \mathfrak{J}_k)\mathcal{L}_{\omega, k}\rangle\nn\\
&\quad + 8 e^{2\delta_1 \lambda^{\frac{1}{3}} t} \mu^{\frac{2}{3}} \sum_{k \in \Z\setminus\{0\}} |k|^{2m-\frac{2}{3}} \operatorname{Re} \langle \partial_{y} \omega_k, (1+ \frac{c_{\beta}}{4}\mathfrak{J}_k) \partial_{y} \mathcal{L}_{\omega, k}\rangle\nn\\
&\quad  +  e^{2\delta_1 \lambda^{\frac{1}{3}} t} \mu^{\frac{1}{3}}  \sum_{k \in \Z\setminus\{0\}} |k|^{2m-\frac{4}{3}} \Big[ \operatorname{Re} \langle ik  \omega_k,  \partial_{y} \mathcal{L}_{\omega, k}\rangle  + \operatorname{Re} \langle ik \mathcal{L}_{\omega, k},  \partial_{y}\omega_k \rangle\Big],
\end{align*}
and
\begin{align*}
\mathcal{N}_{\theta, \neq} &= 32 e^{2\delta_1 \lambda^{\frac{1}{3}} t} \sum_{k \in \Z\setminus\{0\}} |k|^{2m} \operatorname{Re} \langle \theta_k,  \mathcal{N}_{\theta, k}\rangle + 2 e^{2\delta_1 \lambda^{\frac{1}{3}} t} \nu^{\frac{2}{3}} \sum_{k \in \Z\setminus\{0\}} |k|^{2m-\frac{2}{3}} \operatorname{Re} \langle \partial_{y} \theta_k,  \partial_{y} \mathcal{N}_{\theta, k}\rangle\nn\\
&\quad  + e^{2\delta_1 \lambda^{\frac{1}{3}} t}  \nu^{\frac{1}{3}}  \sum_{k \in \Z\setminus\{0\}} |k|^{2m-\frac{4}{3}} \Big[ \operatorname{Re} \langle ik  \theta_k,  \partial_{y} \mathcal{N}_{\theta, k}\rangle  + \operatorname{Re} \langle ik \mathcal{N}_{\theta, k},  \partial_{y}\theta_k \rangle\Big],\nn\\
&:= \mathcal{N}_{\theta, \alpha} + \mathcal{N}_{\theta, \beta} + \mathcal{N}_{\theta, \gamma}.\nn\\
\mathcal{N}_{\omega, \neq} &= 256 e^{2\delta_1 \lambda^{\frac{1}{3}} t} \sum_{k \in \Z\setminus\{0\}} |k|^{2m} \operatorname{Re} \langle \omega_k,  ( 1+ \frac{c_{\alpha}}{128} \mathfrak{J}_k)\mathcal{N}_{\omega, k}\rangle\nn\\
&\quad + 8 e^{2\delta_1 \lambda^{\frac{1}{3}} t} \mu^{\frac{2}{3}} \sum_{k \in \Z\setminus\{0\}} |k|^{2m-\frac{2}{3}} \operatorname{Re} \langle \partial_{y} \omega_k, (1+ \frac{c_{\beta}}{4}\mathfrak{J}_k) \partial_{y} \mathcal{N}_{\omega, k}\rangle\nn\\
&\quad  + e^{2\delta_1 \lambda^{\frac{1}{3}} t} \mu^{\frac{1}{3}}  \sum_{k \in \Z\setminus\{0\}} |k|^{2m-\frac{4}{3}} \Big[ \operatorname{Re} \langle ik  \omega_k,  \partial_{y} \mathcal{N}_{\omega, k}\rangle  + \operatorname{Re} \langle ik \mathcal{N}_{\omega, k},  \partial_{y}\omega_k \rangle\Big]\nn\\
& := \mathcal{N}_{\omega, \alpha} + \mathcal{N}_{\omega, \beta} + \mathcal{N}_{\omega, \gamma}.
\end{align*}
For the linear terms $ \mathcal{L}_{\theta, \neq}$ and $ \mathcal{L}_{\omega, \neq}$, we have already provided estimates in Lemma \ref{le_w6}. Specifically,
$ \mathcal{L}_{\theta, \neq}$ and $ \mathcal{L}_{\omega, \neq}$ satisfy the following estimates:
\begin{align}
	\mathcal{L}_{\theta, \neq} &\le -\frac{1}{4} \mathscr{D}_{\theta, \neq} - \frac{1}{2}\nu^{\frac{1}{3}} e^{2\delta_1 \lambda^{\frac{1}{3}} t} \sum_{k \in \Z\setminus\{0\}} |k|^{2m+\frac{2}{3}} E_{\theta, k}.\label{lta=}\\
	\mathcal{L}_{\omega, \neq} &\le -\frac{1}{16}c_{\beta} \mathscr{D}_{\omega, \neq} - \frac{1}{2}\mu^{\frac{1}{3}}e^{2\delta_1 \lambda^{\frac{1}{3}} t} \sum_{k \in \Z\setminus\{0\}} |k|^{2m+\frac{2}{3}} E_{\omega, k} +C_0 e^{2\delta_1 \lambda^{\frac{1}{3}} t} \mu^{-\frac{1}{3}} \sum_{k \in \Z\setminus\{0\}} |k|^{\frac{4}{3}} \| \theta_k\|^2_{L^2}\nn\\
	&\le -\frac{1}{16}c_{\beta} \mathscr{D}_{\omega, \neq} - \frac{1}{2}\mu^{\frac{1}{3}} e^{2\delta_1 \lambda^{\frac{1}{3}} t} \sum_{k \in \Z\setminus\{0\}} |k|^{2m+\frac{2}{3}} E_{\omega, k} + C_0 \mu^{-\frac{1}{3}} \nu^{-\frac{2}{3}} \mathscr{D}_{\theta, 1}^{\frac{1}{2}}  \mathscr{D}_{\theta, 3}^{\frac{1}{2}}.\label{lw=}
\end{align}
Next, we first estimate the nonlinear terms in the $(\omega, \theta)$ equation. We analyze these three nonlinear terms $ (  \mathcal{N}_{\theta, \alpha}, \mathcal{N}_{\theta, \beta}, \mathcal{N}_{\theta, \gamma})$ separately, leading to the following three specific lemmas
\begin{lemma}\label{le_nonta1}
	Under the conditions of Theorem \ref{thm1.2}, then for any integer $k \neq 0$ and any $t>0$, we have
	\begin{align}\label{Ntaa}
		\mathcal{N}_{\theta, \alpha} \lesssim \nu^{-\frac{1}{2}} \mathscr{D}_{\theta}^{\frac{1}{2}} \mathscr{D}_{\omega}^{\frac{1}{2}} \mathcal{E}_{\theta}^{\frac{1}{2}} + \nu^{-\frac{1}{2}} \mathscr{D}_{\theta}  \mathcal{E}_{\omega}^{\frac{1}{2}} + \nu^{-\frac{1}{6}} \mu^{-\frac{1}{3}} \mathscr{D}_{\theta}^{\frac{1}{2}} \mathscr{D}_{\omega}^{\frac{1}{2}} \mathcal{E}_{\theta}^{\frac{1}{2}}.
	\end{align}
	\begin{proof}
		First, we decompose the $k$-mode of the advection term into the following three parts
		\begin{align}\label{split_na}
			\mathcal{N}_{\theta, \alpha}
			& = -32 e^{2\delta_1 \lambda^{\frac{1}{3}} t}  \sum_{k \in \Z\setminus\{0\}} |k|^{2m} \operatorname{Re} \langle \theta_k, \partial_{y} \psi_{0} \cdot ik \theta_k \rangle + 32 e^{2\delta_1 \lambda^{\frac{1}{3}} t}  \sum_{k \in \Z\setminus\{0\}} |k|^{2m} \operatorname{Re} \langle \theta_k,  ik \psi_{k} \cdot  \partial_{y} \theta_0 \rangle\nn\\
			&\quad -32 e^{2\delta_1 \lambda^{\frac{1}{3}} t}  \sum_{\substack{ k,\ell \in \Z\setminus\{0\} \\ k \neq \ell }} |k|^{2m} \operatorname{Re} \langle \theta_k,  \nabla^T_{k-\ell} \psi_{k-\ell} \cdot \nabla_\ell \theta_\ell \rangle   \nn\\
			& := P_{0, \neq} + P_{\neq, 0} + P_{\neq, \neq}.
		\end{align}
		For the term $ P_{0, \neq}$, we can directly obtain
		\begin{align}\label{p0=}
			 P_{0, \neq}
			  &= 32 e^{2\delta_1 \lambda^{\frac{1}{3}} t}  \sum_{k \in \Z\setminus\{0\}} |k|^{2m} \operatorname{Re} \int_{-1}^{1} ik |\theta_k|^2 \cdot \partial_{y} \psi_{0} \mathrm{d} y  =0.
		\end{align}
		For the term $ P_{\neq, 0}$, applying H\"older's inequality and the Gagliardo–Nirenberg inequality yields
		\begin{align}\label{p=0}
			|P_{\neq, 0} | &\lesssim e^{2\delta_1 \lambda^{\frac{1}{3}} t}  \sum_{k \in \Z\setminus\{0\}} \| |k|^m \theta_k\|_{L^2} \|  |k|^{m+1} \psi_k\|_{L^\infty} \| \partial_{y} \theta_0\|_{L^2}\nn\\
			& \lesssim e^{\delta_1 \lambda^{\frac{1}{3}} t} \nu^{-\frac{1}{2}} \mathscr{D}^{\frac{1}{2}}_{\theta, 0} \mathcal{E}^{\frac{1}{2}}_{\theta}   \Big\{ \sum_{k \in \Z\setminus\{0\}}  \|  |k|^{m+1} \psi_k\|_{L^2} \|  |k|^{m+1} \partial_{y}\psi_k\|_{L^2}\Big\}^{\frac{1}{2}}  \nn\\
			& \lesssim \nu^{-\frac{1}{2}} \mathscr{D}^{\frac{1}{2}}_{\theta, 0} \mathcal{E}^{\frac{1}{2}}_{\theta} \mathscr{D}^{\frac{1}{2}}_{\omega, 4}
		\end{align}
		where we have used the fact
		\begin{align}\label{fact1}
			e^{2\delta_1 \lambda^{\frac{1}{3}} t}  \sum_{k\neq 0} |k|^{2m+2} \|\psi_{k} \|^2_{L^\infty} \lesssim e^{2\delta_1 \lambda^{\frac{1}{3}} t}  \sum_{k \in \Z\setminus\{0\}}  \|  |k|^{m+1} \psi_k\|_{L^2} \|  |k|^{m+1} \partial_{y}\psi_k\|_{L^2} \lesssim \mathscr{D}_{\omega, 4}.
		\end{align}
		For the term $ P_{\neq, \neq}$, we need to decompose the advection term into the following two parts for estimation
		\begin{align}\label{splitp==}
			P_{\neq, \neq}
			& = 32 e^{2\delta_1 \lambda^{\frac{1}{3}} t}  \sum_{\substack{ k,\ell \in \Z\setminus\{0\} \\ k \neq \ell }} |k|^{2m} \operatorname{Re} \langle \theta_k,  \partial_y( i(k-\ell)\psi_{k-\ell} \cdot  \theta_\ell) \rangle \nn\\
&\quad- 32 e^{2\delta_1 \lambda^{\frac{1}{3}} t}  \sum_{\substack{ k,\ell \in \Z\setminus\{0\} \\ k \neq \ell }} |k|^{2m} \operatorname{Re} \langle \theta_k,  ik\partial_y \psi_{k-\ell} \cdot  \theta_\ell \rangle  := 	P_{\neq, \neq, 1} + 	P_{\neq, \neq,2}.
		\end{align}
		For the first term in the above expression, we first apply H\"older's inequality, followed by Young's inequality for convolutions, to obtain
		\begin{align}\label{p==1}
			| P_{\neq, \neq, 1} | &= 32 e^{2\delta_1 \lambda^{\frac{1}{3}} t}  \Big| \sum_{\substack{ k,\ell \in \Z\setminus\{0\} \\ k \neq \ell }} |k|^{2m} \operatorname{Re} \langle \partial_y\theta_k,   i(k-\ell)\psi_{k-\ell} \cdot  \theta_\ell \rangle  \Big|\nn\\
			& \lesssim e^{2\delta_1 \lambda^{\frac{1}{3}} t}  \sum_{\substack{ k,\ell \in \Z\setminus\{0\} \\ k \neq \ell }} |k|^{m}\big( |k-\ell|^m + |\ell|^m\big) \|\partial_{y} \theta_k \|_{L^2} \| (k-\ell) \psi_{k-\ell}\|_{L^\infty} \| \theta_k\|_{L^2}\nn\\
			&\lesssim e^{2\delta_1 \lambda^{\frac{1}{3}} t}  \Big( \sum_{k\neq 0} |k|^{2m} \|\partial_{y} \theta_k \|^2_{L^2}  \Big)^{\frac{1}{2}} \cdot \Big( \sum_{k\neq 0} |k|^{2m+2} \|\psi_{k} \|^2_{L^\infty}  \Big)^{\frac{1}{2}}  \cdot \Big( \sum_{k\neq 0} \|\theta_{k} \|_{L^2}  \Big)\nn\\
			& \quad + e^{2\delta_1 \lambda^{\frac{1}{3}} t}  \Big( \sum_{k\neq 0} |k|^{2m} \|\partial_{y} \theta_k \|^2_{L^2}  \Big)^{\frac{1}{2}} \cdot \Big( \sum_{k\neq 0} |k|^{2m} \|\theta_{k} \|^2_{L^2}  \Big)^{\frac{1}{2}}  \cdot \Big( \sum_{k\neq 0} |k| \|\psi_{k} \|_{L^\infty}  \Big)\nn\\
			&\lesssim \nu^{-\frac{1}{2}} \mathscr{D}_{\theta, 1}^{\frac{1}{2}} \mathscr{D}_{\omega, 4}^{\frac{1}{2}} \cdot \Big( \sum_{k\neq 0} \|\theta_{k} \|_{L^2}  \Big) + \nu^{-\frac{1}{2}} \mathscr{D}_{\theta, 1}^{\frac{1}{2}} \mathcal{E}_{\theta}^{\frac{1}{2}} \cdot \Big( \sum_{k\neq 0} |k|^{2m+2} \|\psi_{k} \|^2_{L^\infty}  \Big)^{\frac{1}{2}}  \Big( \sum_{k\neq 0} |k|^{-2m}   \Big)^{\frac{1}{2}}\nn\\
			& \lesssim  \nu^{-\frac{1}{2}} \mathscr{D}_{\theta, 1}^{\frac{1}{2}} \mathscr{D}_{\omega, 4}^{\frac{1}{2}} \mathcal{E}_{\theta}^{\frac{1}{2}},
 		\end{align}
 		where we have used the facts \eqref{fact1} and $ m >\frac{1}{2}$.

 For the second term $ P_{\neq, \neq,2}$, we need to decompose the frequency $k$. First, let $ A_1(k,\ell):=\{  (k,\ell) \in \Z\setminus\{0\}\big| |k-\ell|\le \frac{|\ell|}{2}\  \text{and} \ k \neq \ell \}$, we define the following two characteristic functions:
 		\begin{align*}
 			\chi_0(k,\ell)= \begin{cases}1, & \text { for } (k,\ell) \in A_1 \\ 0, & \text { for } (k,\ell) \notin A_1\end{cases} \quad \text{and} \quad
 			\chi_1(k,\ell)= \begin{cases}0, & \text { for } (k,\ell) \in A_1 \\ 1, & \text { for } (k,\ell) \notin A_1\end{cases}
 		\end{align*}
 		Here, we can decompose the term $ P_{\neq, \neq,2}$ into the following two parts
		\begin{align}\label{p==2}
			|P_{\neq, \neq,2} | &\lesssim e^{2\delta_1 \lambda^{\frac{1}{3}} t}  \sum_{\substack{ k,\ell \in \Z\setminus\{0\} \\ k \neq \ell }} |k|^{2m+1} \| \theta_{k}\|_{L^2} \| \partial_{y} \psi_{k-\ell}\|_{L^\infty} \|\theta_{\ell} \|_{L^2}\nn\\
			&\lesssim e^{2\delta_1 \lambda^{\frac{1}{3}} t}  \sum_{\substack{ k,\ell \in \Z\setminus\{0\} \\ k \neq \ell }} \Big( \chi_0(k,\ell) + \chi_1(k,\ell) \Big)|k|^{2m+1} \| \theta_{k}\|_{L^2} \| \partial_{y} \psi_{k-\ell}\|_{L^\infty} \|\theta_{\ell} \|_{L^2}\nn\\
			&:= P_{\neq, \neq, \chi_0} + P_{\neq, \neq, \chi_1}.
		\end{align}
		For $ (k,\ell)\in A_1$, we have $ |k|\lesssim |\ell|$. Directly applying H\"older's inequality and Young's inequality, we obtain
		\begin{align}\label{p==x0}
			| P_{\neq, \neq, \chi_0}| &\lesssim e^{2\delta_1 \lambda^{\frac{1}{3}} t}  \sum_{\substack{ k,\ell \in \Z\setminus\{0\} \\ k \neq \ell }} |k|^{m+\frac{1}{3}} \|\theta_{k} \|_{L^2} |\ell|^{m+\frac{2}{3}} \|\theta_{\ell} \|_{L^2} \| \partial_{y} \psi_{k-\ell}\|_{L^\infty} \nn\\
			&\lesssim e^{\delta_1 \lambda^{\frac{1}{3}} t} \nu^{-\frac{1}{6}} \mathscr{D}_{\theta,3}^{\frac{1}{2}} \Big( \sum_{k\neq 0} |k|^{2m+\frac{4}{3}} \| \theta_{k}\|^2_{L^2}\Big)^{\frac{1}{2}} \Big( \sum_{k\neq 0} \| \partial_{y} \psi_{k}\|_{L^\infty}\Big)\nn\\
			& \lesssim \nu^{-\frac{1}{2}} \mathscr{D}_{\theta,3}^{\frac{3}{4}} \mathscr{D}_{\theta,1}^{\frac{1}{4}} \mathcal{E}^{\frac{1}{2}}_{\omega},
		\end{align}
		where we have used the facts
		\begin{align}\label{fact3}
			e^{2\delta_1 \lambda^{\frac{1}{3}} t}  \sum_{k\neq 0} |k|^{2m+\frac{4}{3}} \| \theta_{k}\|^2_{L^2} &\lesssim e^{2\delta_1 \lambda^{\frac{1}{3}} t}  \nu^{-\frac{2}{3}}\sum_{k\neq 0} \Big( \nu^{\frac{1}{6}}|k|^{m+\frac{1}{3}} \| \theta_{k}\|_{L^2}\Big) \Big( \nu^{\frac{1}{2}}|k|^{m} \| \nabla_k \theta_{k}\|_{L^2}\Big)\nn\\
			& \lesssim \nu^{-\frac{2}{3}} \mathscr{D}_{\theta,3}^{\frac{1}{2}} \mathscr{D}_{\theta,1}^{\frac{1}{2}},
		\end{align}
		and
		\begin{align*}
			e^{\delta_1 \lambda^{\frac{1}{3}} t} \sum_{k\neq 0} \| \partial_{y} \psi_{k}\|_{L^\infty} &\lesssim e^{\delta_1 \lambda^{\frac{1}{3}} t} 	\sum_{k\neq 0} |k|^{-m} \Big( |k|^{\frac{m}{2}} \| \partial_{y} \psi_{k}\|^{\frac{1}{2}}_{L^2} \Big) \Big( |k|^{\frac{m}{2}}  \| \partial_{yy} \psi_{k}\|^{\frac{1}{2}}_{L^2} \Big) \nn\\
			&\lesssim e^{\delta_1 \lambda^{\frac{1}{3}} t} \sum_{k\neq 0} |k|^{-m} \Big( |k|^{m} \| \partial_{yy} \psi_{k}\|_{L^2} \Big) \lesssim \mathcal{E}^{\frac{1}{2}}_{\omega}.
		\end{align*}
			For $ (k,\ell)\notin A_1$, we have $ |k|\lesssim |k-\ell|+ |\ell|\lesssim |k-\ell|$. Directly applying H\"older's inequality and Young's inequality, we obtain
		\begin{align}\label{p==x1}
			| P_{\neq, \neq, \chi_1}|  &\lesssim e^{2\delta_1 \lambda^{\frac{1}{3}} t}  \sum_{\substack{ k,\ell \in \Z\setminus\{0\} \\ k \neq \ell }} |k|^{m+\frac{1}{3}} \|\theta_{k} \|_{L^2} |k-\ell|^{m+\frac{2}{3}} \| \partial_{y} \psi_{k-\ell}\|_{L^\infty}  \|\theta_{\ell} \|_{L^2} \nn\\
			&\lesssim e^{\delta_1 \lambda^{\frac{1}{3}} t} \nu^{-\frac{1}{6}} \mathscr{D}_{\theta,3}^{\frac{1}{2}} \Big( \sum_{k\neq 0} |k|^{2m+\frac{4}{3}} \| \partial_{y} \psi_{k}\|^2_{L^\infty} \Big)^{\frac{1}{2}} \Big( \sum_{k\neq 0} \| \theta_{k}\|_{L^2}\Big)\nn\\
			& \lesssim \nu^{-\frac{1}{6}} \mu^{-\frac{1}{3}} \mathscr{D}_{\theta,3}^{\frac{1}{2}} \mathscr{D}_{\omega,5}^{\frac{1}{2}} \mathcal{E}^{\frac{1}{2}}_{\theta},
		\end{align}
		where
		\begin{align}\label{fact2}
			e^{2\delta_1 \lambda^{\frac{1}{3}} t} \sum_{k\neq 0} |k|^{2m+\frac{7}{3}} \| \partial_{y} \psi_{k}\|^2_{L^\infty} \lesssim e^{2\delta_1 \lambda^{\frac{1}{3}} t} \sum_{k\neq 0} |k|^{m+\frac{2}{3}} \| |k|\partial_{y} \psi_{k}\|_{L^2}\cdot  |k|^{m+\frac{2}{3}}\| \partial_{yy} \psi_{k}\|_{L^2} \lesssim \mu^{-\frac{2}{3}}\mathscr{D}_{\omega, 5}.
		\end{align}
		Collecting \eqref{split_na}--\eqref{p==x1}, we complete the proof of the lemma.
	\end{proof}
\end{lemma}

\begin{lemma}\label{le_nonta2}
	Under the conditions of Theorem \ref{thm1.2}, then for any integer $k \neq 0$ and any $t>0$, we have
	\begin{align}
		\mathcal{N}_{\theta, \beta}  \lesssim \nu^{-\frac{1}{2}} \mathscr{D}_{\theta}^{\frac{1}{2}} \mathscr{D}_{\omega}^{\frac{1}{2}} \mathcal{E}_{\theta}^{\frac{1}{2}} + (\nu^{-\frac{1}{2}} + \nu^{-\frac{1}{3}} +\mu^{-\frac{1}{6}}) \mathscr{D}_{\theta}  \mathcal{E}_{\omega}^{\frac{1}{2}} + \nu^{-\frac{1}{6}} \mu^{-\frac{1}{3}} \mathscr{D}_{\theta}^{\frac{1}{2}} \mathscr{D}_{\omega}^{\frac{1}{2}} \mathcal{E}_{\theta}^{\frac{1}{2}}.
	\end{align}
	\begin{proof}
		First, following a method similar to \eqref{split_na}, we decompose the frequency $(k,\ell)$ into three cases: $k=\ell$, $\ell=0$, and $k \neq \ell$ with $\ell \neq 0$. Specifically, we have
		\begin{align}\label{ntab}
			\mathcal{N}_{\theta, \beta} &= -2 e^{2\delta_1 \lambda^{\frac{1}{3}} t}  \nu^{\frac{2}{3}} \sum_{k \in \Z\setminus\{0\}} |k|^{2m-\frac{2}{3}} \operatorname{Re} \langle \partial_{y} \theta_k, \partial_{y}\Big[ \sum_{\ell \in \Z} \partial_{y} \psi_{k-\ell} \cdot i\ell \theta_\ell\Big] \rangle \nn\\
			&\quad + 2 e^{2\delta_1 \lambda^{\frac{1}{3}} t}  \nu^{\frac{2}{3}} \sum_{k \in \Z\setminus\{0\}} |k|^{2m-\frac{2}{3}} \operatorname{Re} \langle \partial_{y} \theta_k, \partial_{y} \Big[ \sum_{\ell \in \Z} i(k-\ell) \psi_{k-\ell} \cdot \partial_{y} \theta_\ell \Big]\rangle\nn\\
			& = -2 e^{2\delta_1 \lambda^{\frac{1}{3}} t}  \nu^{\frac{2}{3}} \sum_{k \in \Z\setminus\{0\}} |k|^{2m-\frac{2}{3}} \Bigg(\operatorname{Re} \langle \partial_{y} \theta_k, \partial_{y}\big[ \partial_{y} \psi_{0} \cdot ik \theta_k\big] \rangle - \operatorname{Re} \langle \partial_{y} \theta_k,  \partial_{y} \big[ ik \psi_{k} \cdot  \partial_{y} \theta_0\big] \rangle \Bigg)   \nn\\
			&\quad -2 e^{2\delta_1 \lambda^{\frac{1}{3}} t}  \nu^{\frac{2}{3}} \sum_{\substack{ k,\ell \in \Z\setminus\{0\} \\ k \neq \ell }} |k|^{2m-\frac{2}{3}} \operatorname{Re} \langle \partial_{y} \theta_k,  \partial_{y} \big[ \nabla^\top_{k-\ell} \psi_{k-\ell} \cdot \nabla_\ell \theta_\ell\big] \rangle   \nn\\
			& := Q_{0, \neq} + Q_{\neq, 0} + Q_{\neq, \neq}.
		\end{align}
		For the case $k=\ell$, we decompose it into the following two terms
		\begin{align}\label{split_qo=}
			Q_{0, \neq} &= -2 e^{2\delta_1 \lambda^{\frac{1}{3}} t} \nu^{\frac{2}{3}} \sum_{k \in \Z\setminus\{0\}} |k|^{2m-\frac{2}{3}} \Big(\operatorname{Re} \langle \partial_{y} \theta_k,  \partial_{y} \psi_{0} \cdot ik \partial_{y}\theta_k \rangle + \operatorname{Re} \langle \partial_{y} \theta_k,  \partial_{yy} \psi_{0} \cdot ik \theta_k \rangle \Big)   \nn\\
			&:= Q_{0, \neq, 1} + Q_{0, \neq, 2}.
		\end{align}
		For the first term in the above expression, we can directly obtain $ Q_{0, \neq, 1} =0$. For the second term, applying H\"older's inequality and the Gagliardo–Nirenberg inequality, we obtain
		\begin{align}\label{q0=2}
			|Q_{0, \neq, 2}| &\lesssim  \sum_{k \in \Z\setminus\{0\}}  \nu^{\frac{1}{2}} |k|^{m} \|\partial_{y} \theta_{k} \|_{L^2} \cdot \nu^{\frac{1}{6}}|k|^{m+\frac{1}{3}}\| \theta_{k}\|_{L^2} \|\partial_{yy}\psi_{0} \|_{L^\infty}\nn\\
			&\lesssim \mathscr{D}_{\theta, 1}^{\frac{1}{2}} \mathscr{D}_{\theta, 3}^{\frac{1}{2}} \| \omega_0\|_{L^\infty}\lesssim \mathscr{D}_{\theta, 1}^{\frac{1}{2}} \mathscr{D}_{\theta, 3}^{\frac{1}{2}} \| \omega_0\|_{L^2}^{\frac{1}{2}} \|\partial_{y} \omega_0\|_{L^2}^{\frac{1}{2}}\nn\\
			&\lesssim \mu^{-\frac{1}{6}}\mathscr{D}_{\theta, 1}^{\frac{1}{2}} \mathscr{D}_{\theta, 3}^{\frac{1}{2}} \mathcal{E}_{\omega, 0}^{\frac{1}{2}}.
		\end{align}
		For the case $\ell=0$, applying H\"older's inequality and \eqref{fact1} yields
		\begin{align}\label{q=0}
			| Q_{\neq, 0}| &= 2 e^{2\delta_1 \lambda^{\frac{1}{3}} t} \Big| \nu^{\frac{2}{3}} \sum_{k \in \Z\setminus\{0\}} |k|^{2m-\frac{2}{3}} \operatorname{Re} \langle \partial_{yy} \theta_k,   ik \psi_{k} \cdot  \partial_{y} \theta_0 \rangle\Big|\nn\\
			&\lesssim e^{2\delta_1 \lambda^{\frac{1}{3}} t} \nu^{-\frac{1}{6}} \sum_{k \in \Z\setminus\{0\}} \nu^{\frac{5}{6}}|k|^{m-\frac{1}{3}} \| \partial_{yy} \theta_k\|_{L^2} \cdot |k|^{m+\frac{2}{3}} \| \psi_{k}\|_{L^\infty} \| \partial_{y} \theta_0\|_{L^2}\nn\\
			&\lesssim \nu^{-\frac{1}{2}} \mathscr{D}_{\theta, 2}^{\frac{1}{2}} \mathscr{D}_{\omega, 4}^{\frac{1}{2}} \mathcal{E}_{\theta, 0}^{\frac{1}{2}}.
		\end{align}
		For the case $k \neq \ell$ and $\ell\neq 0$, we handle the first term using integration by parts and expand the second term, obtaining the following three terms
		\begin{align}\label{split_q==}
			Q_{\neq, \neq} &= -2 e^{2\delta_1 \lambda^{\frac{1}{3}} t} \nu^{\frac{2}{3}} \sum_{\substack{ k,\ell \in \Z\setminus\{0\} \\ k \neq \ell }} |k|^{2m-\frac{2}{3}} \operatorname{Re}  \Big\langle \partial_{y} \theta_k,  \partial_{y} \big( \partial_{y} \psi_{k-\ell} \cdot i\ell \theta_\ell\big) -\partial_{y} \big(  i(k-\ell)\psi_{k-\ell} \cdot  \partial_{y}\theta_\ell\big) \Big\rangle\nn\\
			&=2 e^{2\delta_1 \lambda^{\frac{1}{3}} t} \nu^{\frac{2}{3}} \sum_{\substack{ k,\ell \in \Z\setminus\{0\} \\ k \neq \ell }} |k|^{2m-\frac{2}{3}} \Big[ \operatorname{Re} \langle \partial_{yy} \theta_k,   \partial_{y} \psi_{k-\ell} \cdot i\ell \theta_\ell \rangle +  \operatorname{Re} \langle \partial_{y} \theta_k,   i(k-\ell)\partial_{y}\psi_{k-\ell} \cdot  \partial_{y}\theta_\ell\rangle\Big] \nn\\
			& \quad + 2 e^{2\delta_1 \lambda^{\frac{1}{3}} t} \nu^{\frac{2}{3}} \sum_{\substack{ k,\ell \in \Z\setminus\{0\} \\ k \neq \ell }} |k|^{2m-\frac{2}{3}} \operatorname{Re} \langle \partial_{y} \theta_k,   i(k-\ell)\psi_{k-\ell} \cdot  \partial_{yy}\theta_\ell\rangle\nn\\
			& := Q_{\neq, \neq, 1} + Q_{\neq, \neq,2} + Q_{\neq, \neq,3}.
		\end{align}
		For the second term $ Q_{\neq, \neq,2}$, we use the two characteristic functions defined in the Lemma \ref{le_nonta1} to decompose it into the following two parts
		\begin{align}\label{q==2}
			Q_{\neq, \neq,2} &=  2 e^{2\delta_1 \lambda^{\frac{1}{3}} t} \nu^{\frac{2}{3}} \sum_{\substack{ k,\ell \in \Z\setminus\{0\} \\ k \neq \ell }} \Big( \chi_0(k,\ell) + \chi_1(k,\ell) \Big)|k|^{2m-\frac{2}{3}} \operatorname{Re} \langle \partial_{y} \theta_k,   i(k-\ell)\partial_{y}\psi_{k-\ell} \cdot  \partial_{y}\theta_\ell\rangle\nn\\
			&:= Q_{2, \chi_0} + Q_{2, \chi_1}.
		\end{align}
			For $ (k,\ell)\in A_1$, we have $ |k|\lesssim |\ell|$. Directly applying  Young's inequality and \eqref{fact2}, we obtain
		\begin{align}\label{q2x0}
			Q_{2, \chi_0} &\lesssim  e^{2\delta_1 \lambda^{\frac{1}{3}} t} \nu^{-\frac{1}{6}}\sum_{\substack{ k,\ell \in \Z\setminus\{0\} \\ k \neq \ell }} \nu^{\frac{1}{3}} |k|^{m-\frac{1}{3}} \| \partial_{y} \theta_{k}\|_{L^2} \nu^{\frac{1}{2}} |\ell|^{m-\frac{1}{3}} \| \partial_{y} \theta_{\ell}\|_{L^2} \||k-\ell|\partial_{y} \psi_{k-\ell} \|_{L^\infty}\nn\\
			&\lesssim \nu^{-\frac{1}{6}} \mathscr{D}_{\theta, 1}^{\frac{1}{2}} \mathcal{E}_{\theta}^{\frac{1}{2}} \sum_{k\neq 0} |k|^{-m-\frac{1}{6}}|k|^{m+\frac{7}{6}} \|\partial_{y} \psi_{k} \|_{L^\infty}\lesssim \nu^{-\frac{1}{6}} \mu^{-\frac{1}{3}} \mathscr{D}_{\theta, 1}^{\frac{1}{2}} \mathcal{E}_{\theta}^{\frac{1}{2}} \mathscr{D}_{\omega, 5}^{\frac{1}{2}}.
 		\end{align}
 			For $ (k,\ell)\notin A_1$, we have $ |k|\lesssim |k-\ell|+ |\ell|\lesssim |k-\ell|$. Similarly, we can get
 		\begin{align}\label{q2x1}
 			Q_{2, \chi_1} &\lesssim e^{2\delta_1 \lambda^{\frac{1}{3}} t} \nu^{\frac{1}{3}} \sum_{\substack{ k,\ell \in \Z\setminus\{0\} \\ k \neq \ell }} \nu^{\frac{1}{3}} |k|^{m-\frac{1}{3}} \| \partial_{y} \theta_{k}\|_{L^2}  \| \partial_{y} \theta_{\ell}\|_{L^2}  |k-\ell|^{m+\frac{2}{3}}\|\partial_{y} \psi_{k-\ell} \|_{L^\infty}\nn\\
 			&\lesssim \nu^{-\frac{1}{6}} \mu^{-\frac{1}{3}} \mathscr{D}_{\theta, 1}^{\frac{1}{2}} \mathcal{E}_{\theta}^{\frac{1}{2}} \mathscr{D}_{\omega, 5}^{\frac{1}{2}}.
 		\end{align}
 		For the term $ Q_{\neq, \neq,3} $, we handle it similarly to $ Q_{\neq, \neq,2} $ by decomposing it into two parts
 		\begin{align*}
 			Q_{\neq, \neq,3} &=  2 e^{2\delta_1 \lambda^{\frac{1}{3}} t} \nu^{\frac{2}{3}} \sum_{\substack{ k,\ell \in \Z\setminus\{0\} \\ k \neq \ell }} \Big( \chi_0(k,\ell) + \chi_1(k,\ell) \Big)|k|^{2m-\frac{2}{3}} \operatorname{Re} \langle \partial_{y} \theta_k,   i(k-\ell)\psi_{k-\ell} \cdot  \partial_{yy}\theta_\ell\rangle\nn\\
 			&:= Q_{3, \chi_0} + Q_{3, \chi_1}.
 		\end{align*}
 		For the first term, similar to the estimation of $ Q_{2, \chi_0}$ and using \eqref{fact1}, we obtain
 		\begin{align}\label{q3x0}
 			Q_{3, \chi_0} &\lesssim e^{2\delta_1 \lambda^{\frac{1}{3}} t} \nu^{-\frac{1}{2}} \sum_{\substack{ k,\ell \in \Z\setminus\{0\} \\ k \neq \ell }} \nu^{\frac{1}{3}} |k|^{m-\frac{1}{3}} \| \partial_{y} \theta_{k}\|_{L^2} |k-\ell| \| \psi_{k-\ell}\|_{L^\infty} \cdot \nu^{\frac{5}{6}} |\ell|^{m-\frac{1}{3}} \| \partial_{yy} \theta_{\ell}\|_{L^2}\nn\\
 			&\lesssim \nu^{-\frac{1}{2}} \mathscr{D}_{\theta, 2}^{\frac{1}{2}} \mathcal{E}_{\theta}^{\frac{1}{2}} \mathscr{D}_{\omega, 4}^{\frac{1}{2}}.
 		\end{align}
 		For the second term, similar to the estimation of $ Q_{2, \chi_1}$ and using \eqref{fact1}, we can get
 		\begin{align}\label{q3x1}
 			Q_{3, \chi_1} &\lesssim e^{2\delta_1 \lambda^{\frac{1}{3}} t} \nu^{\frac{1}{3}} \sum_{\substack{ k,\ell \in \Z\setminus\{0\} \\ k \neq \ell }} \nu^{\frac{1}{3}} |k|^{m-\frac{1}{3}} \| \partial_{y} \theta_{k}\|_{L^2} |k-\ell|^{m+\frac{2}{3}} \| \psi_{k-\ell}\|_{L^\infty}   \| \partial_{yy} \theta_{\ell}\|_{L^2}\nn\\
 			&\lesssim e^{\delta_1 \lambda^{\frac{1}{3}} t}\nu^{\frac{1}{3}} \mathcal{E}_{\theta}^{\frac{1}{2}} \Big( \sum_{k\neq 0} \| \partial_{yy} \theta_{k}\|^2_{L^2}\Big)^{\frac{1}{2}} \sum_{k\neq 0} |k|^{m+\frac{2}{3}} \| \psi_{k}\|_{L^\infty} \nn\\
 			&\lesssim \nu^{-\frac{1}{2}} \mathcal{E}_{\theta}^{\frac{1}{2}} \mathscr{D}_{\theta,2}^{\frac{1}{2}} \sum_{k\neq 0} |k|^{m+\frac{1}{6}} \| \nabla_k \psi_{k}\|_{L^2}\nn\\
 			&\lesssim \nu^{-\frac{1}{2}} \mathscr{D}_{\theta, 2}^{\frac{1}{2}} \mathcal{E}_{\theta}^{\frac{1}{2}} \mathscr{D}_{\omega, 4}^{\frac{1}{2}}.
 		\end{align}
 		For the term $ Q_{\neq, \neq,1} $, we handle it similarly to $ Q_{\neq, \neq,2} $ by decomposing it into two parts
 		\begin{align*}
 			Q_{\neq, \neq,1} &=  2 e^{2\delta_1 \lambda^{\frac{1}{3}} t} \nu^{\frac{2}{3}} \sum_{\substack{ k,\ell \in \Z\setminus\{0\} \\ k \neq \ell }} \Big( \chi_0(k,\ell) + \chi_1(k,\ell) \Big)|k|^{2m-\frac{2}{3}} \operatorname{Re} \langle \partial_{yy} \theta_k,   \partial_{y} \psi_{k-\ell} \cdot  i\ell \theta_\ell\rangle\nn\\
 			&:= Q_{1, \chi_0} + Q_{1, \chi_1}.
 		\end{align*}
 		For the first term above, using Young's inequality along with \eqref{fact3} and \eqref{w_fact}, we obtain
 		\begin{align}\label{q1x0}
 			Q_{1, \chi_0} &\lesssim e^{2\delta_1 \lambda^{\frac{1}{3}} t} \nu^{-\frac{1}{6}} \sum_{\substack{ k,\ell \in \Z\setminus\{0\} \\ k \neq \ell }} \nu^{\frac{5}{6}}|k|^{m-\frac{1}{3}} \| \partial_{yy} \theta_k\|_{L^2} |\ell|^{m+\frac{2}{3}} \|  \theta_\ell\|_{L^2} \| \partial_{y} \psi_{k-\ell}\|_{L^\infty}\nn\\
 			&\lesssim e^{\delta_1 \lambda^{\frac{1}{3}} t} \nu^{-\frac{1}{6}} \mathscr{D}_{\theta, 2}^{\frac{1}{2}}  \Big( \sum_{k\neq 0} |k|^{2m+\frac{4}{3}} \|\theta_{k} \|^2_{L^2}\Big)^{\frac{1}{2}} \sum_{k\neq 0} \| \partial_{y} \psi_{k}\|_{L^\infty}\nn\\
 			&\lesssim \nu^{-\frac{1}{2}} \mathscr{D}_{\theta, 2}^{\frac{1}{2}} \mathscr{D}_{\theta, 1}^{\frac{1}{4}} \mathscr{D}_{\theta, 3}^{\frac{1}{4}} \sum_{k\neq 0} \| \omega_k\|_{L^2} \lesssim \nu^{-\frac{1}{2}} \mathscr{D}_{\theta, 2}^{\frac{1}{2}} \mathscr{D}_{\theta, 1}^{\frac{1}{4}} \mathscr{D}_{\theta, 3}^{\frac{1}{4}} \mathcal{E}_{\omega}^{\frac{1}{2}}.
 		\end{align}
 		Similarly, for the second term, we also obtain from \eqref{fact3} that
 		\begin{align}\label{q1x1}
 			Q_{1, \chi_1}  &\lesssim e^{2\delta_1 \lambda^{\frac{1}{3}} t} \nu^{-\frac{1}{6}} \sum_{\substack{ k,\ell \in \Z\setminus\{0\} \\ k \neq \ell }} \nu^{\frac{5}{6}}|k|^{m-\frac{1}{3}} \| \partial_{yy} \theta_k\|_{L^2} |k-\ell|^{m-\frac{1}{3}}  \| \partial_{y} \psi_{k-\ell}\|_{L^\infty} \| |\ell| \theta_\ell\|_{L^2}\nn\\
 			&\lesssim e^{\delta_1 \lambda^{\frac{1}{3}} t} \nu^{-\frac{1}{6}} \mathscr{D}_{\theta, 2}^{\frac{1}{2}}  \Big( \sum_{k\neq 0} |k|^{2} \|\theta_{k} \|^2_{L^2}\Big)^{\frac{1}{2}} \sum_{k\neq 0} |k|^{m-\frac{1}{3}}\| \partial_{y} \psi_{k}\|_{L^\infty}\nn\\
 			&\lesssim \nu^{-\frac{1}{3}} \mathscr{D}_{\theta, 2}^{\frac{1}{2}} \mathscr{D}_{\theta, 3}^{\frac{1}{2}}  \mathcal{E}_{\omega}^{\frac{1}{2}}.
 		\end{align}
 		In the last inequality, we used the fact that
 		\begin{align*}
 			e^{\delta_1 \lambda^{\frac{1}{3}} t}\sum_{k\neq 0} |k|^{m-\frac{1}{3}}\| \partial_{y} \psi_{k}\|_{L^\infty} \lesssim e^{\delta_1 \lambda^{\frac{1}{3}} t}\sum_{k\neq 0} |k|^{m-\frac{1}{3}} \| \omega_k\|^{\frac{1}{2}}_{L^2}  \| \partial_{y} \psi_{k}\|^{\frac{1}{2}}_{L^2} \lesssim e^{\delta_1 \lambda^{\frac{1}{3}} t}\sum_{k\neq 0} |k|^{m-\frac{5}{6}} \| \omega_k\|_{L^2} \lesssim   \mathcal{E}_{\omega}^{\frac{1}{2}}.
 		\end{align*}
 		Collecting \eqref{ntab}-\eqref{q1x1}, we complete the proof of the lemma.
	\end{proof}
\end{lemma}

\begin{lemma}\label{le_nonta3}
	Under the conditions of Theorem \ref{thm1.2}, then for any integer $k \neq 0$ and any $t>0$, we have
	\begin{align}
	\mathcal{N}_{\theta, \gamma}  \lesssim \nu^{-\frac{1}{2}} \mathscr{D}_{\theta}^{\frac{1}{2}} \mathscr{D}_{\omega}^{\frac{1}{2}} \mathcal{E}_{\theta}^{\frac{1}{2}} + ( \nu^{-\frac{1}{3}} +\mu^{-\frac{1}{6}}) \mathscr{D}_{\theta}  \mathcal{E}_{\omega}^{\frac{1}{2}} + (\nu^{-\frac{1}{6}} \mu^{-\frac{1}{6}} + \mu^{-\frac{1}{6}}) \mathscr{D}_{\theta}^{\frac{1}{2}} \mathscr{D}_{\omega}^{\frac{1}{2}} \mathcal{E}_{\theta}^{\frac{1}{2}}.
	\end{align}
	\begin{proof}
		First, following a method similar to \eqref{split_na}, we decompose the frequency $(k,\ell)$ into three cases: $k=\ell$, $\ell=0$, and $k \neq \ell$ with $\ell \neq 0$. Specifically, we have
		\begin{align}
			\mathcal{N}_{\theta, \gamma} &= -2 e^{2\delta_1 \lambda^{\frac{1}{3}} t} \nu^{\frac{1}{3}} \sum_{k\neq 0} |k|^{2m - \frac{4}{3}}\operatorname{Re} \langle ik \theta_{k}, \partial_{y}(\nabla^\top \psi \cdot \nabla \theta)_k \rangle\nn\\
			& = -2 e^{2\delta_1 \lambda^{\frac{1}{3}} t} \nu^{\frac{1}{3}} \sum_{k\neq 0} |k|^{2m - \frac{4}{3}}\operatorname{Re} \Big\langle ik \theta_{k}, \partial_{y}\big( \partial_{y} \psi_{0} \cdot ik \theta_{k}\big) \Big\rangle \nn\\
&\quad+ 2 e^{2\delta_1 \lambda^{\frac{1}{3}} t} \nu^{\frac{1}{3}} \sum_{k\neq 0} |k|^{2m - \frac{4}{3}}\operatorname{Re} \Big\langle ik \theta_{k}, \partial_{y}\big(  ik\psi_{k} \cdot  \partial_{y}\theta_{0}\big) \Big\rangle\nn\\
			&\quad -2 e^{2\delta_1 \lambda^{\frac{1}{3}} t} \nu^{\frac{1}{3}} \sum_{\substack{ k,\ell \in \Z\setminus\{0\} \\ k \neq \ell }} |k|^{2m - \frac{4}{3}} \operatorname{Re} \Big\langle ik\theta_k,  \partial_{y}\big(\nabla^\top_{k-\ell} \psi_{k-\ell} \cdot \nabla_\ell \theta_\ell \big) \Big\rangle\nn\\
			& := R_{0, \neq} + R_{\neq, 0} + R_{\neq, \neq}.
		\end{align}
		For the case $k=\ell$, in the above expression, we first expand this term and then apply integration by parts to obtain
		\begin{align}\label{r0=}
			R_{0, \neq} &= -2 e^{2\delta_1 \lambda^{\frac{1}{3}} t} \nu^{\frac{1}{3}} \sum_{k\neq 0} |k|^{2m - \frac{4}{3}}\operatorname{Re} \langle ik \theta_{k},  \partial_{yy} \psi_{0} \cdot ik \theta_{k} \rangle \nn\\
&\quad- 2 e^{2\delta_1 \lambda^{\frac{1}{3}} t} \nu^{\frac{1}{3}} \sum_{k\neq 0} |k|^{2m - \frac{4}{3}}\operatorname{Re} \langle ik \theta_{k},  \partial_{y} \psi_{0} \cdot ik \partial_{y}\theta_{k} \rangle\nn\\
			& = - e^{2\delta_1 \lambda^{\frac{1}{3}} t} \nu^{\frac{1}{3}} \sum_{k\neq 0} |k|^{2m + \frac{2}{3}}\operatorname{Re} \langle  \theta_{k},  \partial_{yy} \psi_{0} \cdot \theta_{k} \rangle \lesssim \mathscr{D}_{\theta, 3} \| \partial_{yy} \psi_{0}\|_{L^\infty} \nn\\
			&\lesssim \mathscr{D}_{\theta, 3} \| \omega_0\|^{\frac{1}{2}}_{L^2} \|\partial_{y}\omega_0 \|^{\frac{1}{2}}_{L^2} \lesssim \mu^{-\frac{1}{6}} \mathscr{D}_{\theta, 3} \mathcal{E}_{\omega}^{\frac{1}{2}}.
		\end{align}
			For the case $\ell=0$, we first apply integration by parts and then use H\"older's inequality along with  \eqref{fact1} to obtain
		\begin{align}\label{r=0}
			R_{0, \neq} &= 2e^{2\delta_1 \lambda^{\frac{1}{3}} t} \nu^{\frac{1}{3}} \sum_{k\neq 0} |k|^{2m - \frac{4}{3}}\operatorname{Re} \langle ik \partial_{y} \theta_{k},  ik\psi_{k} \cdot  \partial_{y}\theta_{0} \rangle \nn\\
			&\lesssim e^{2\delta_1 \lambda^{\frac{1}{3}} t}  \nu^{-\frac{1}{2}} \sum_{k\neq 0}  \nu^{\frac{1}{2}} |k|^{m}\| \partial_{y} \theta_{k}\|_{L^2}  |k|^{m + \frac{2}{3}}\| \psi_{k}\|_{L^\infty} \cdot \nu^{\frac{1}{3}}\| \partial_{y}\theta_{0}\|_{L^2}\nn\\
			&\lesssim \nu^{-\frac{1}{2}} \mathscr{D}_{\theta, 1}^{\frac{1}{2}} \mathscr{D}_{\omega, 4}^{\frac{1}{2}} \mathcal{E}_{\theta}^{\frac{1}{2}}.
		\end{align}
		For the case $k \neq \ell$ and $\ell \neq 0$, we expand the advection term into the following two parts
		\begin{align}\label{split_r==}
				R_{\neq, \neq} &= -2 e^{2\delta_1 \lambda^{\frac{1}{3}} t} \nu^{\frac{1}{3}} \sum_{\substack{ k,\ell \in \Z\setminus\{0\} \\ k \neq \ell }} |k|^{2m - \frac{4}{3}} \operatorname{Re} \Big\langle ik\theta_k,  \partial_{y}\big(\partial_{y} \psi_{k-\ell} \cdot i\ell \theta_\ell \big) \Big\rangle \nn\\
				&\quad + 2 e^{2\delta_1 \lambda^{\frac{1}{3}} t} \nu^{\frac{1}{3}} \sum_{\substack{ k,\ell \in \Z\setminus\{0\} \\ k \neq \ell }} |k|^{2m - \frac{4}{3}} \operatorname{Re} \Big\langle ik\theta_k,  \partial_{y}\big( i(k-\ell)\psi_{k-\ell} \cdot  \partial_{y} \theta_\ell \big) \Big\rangle\nn\\
 &:= R_{\neq, \neq, 1} + R_{\neq, \neq, 2}.
		\end{align}
		For the term $ R_{\neq, \neq, 2}$, applying integration by parts and using \eqref{fact1} yields
		\begin{align}\label{r==2}
			R_{\neq, \neq, 2} &\lesssim e^{2\delta_1 \lambda^{\frac{1}{3}} t} \nu^{\frac{1}{3}} \sum_{\substack{ k,\ell \in \Z\setminus\{0\} \\ k \neq \ell }} |k|^{m - \frac{1}{3}} \| \partial_{y} \theta_{k}\|_{L^2} \big( |k-\ell|^m + |\ell|^m\big)\cdot |k-\ell|\| \psi_{k-\ell}\|_{L^\infty} \| \partial_{y} \theta_{\ell}\|_{L^2}\nn\\
			&\lesssim e^{2\delta_1 \lambda^{\frac{1}{3}} t}\Big(  \sum_{k\neq 0} \nu^{\frac{2}{3}} |k|^{2m - \frac{2}{3}} \| \partial_{y} \theta_{k}\|^2_{L^2}  \Big)^{\frac{1}{2}} \Big( \sum_{k\neq 0}  |k|^{2m +2 } \| \psi_{k}\|^2_{L^\infty} \Big)^{\frac{1}{2}} \Big( \sum_{k\neq 0} \| \partial_{y} \theta_{k}\|_{L^2} \Big)\nn\\
			&\quad + e^{2\delta_1 \lambda^{\frac{1}{3}} t} \Big(  \sum_{k\neq 0} \nu^{\frac{2}{3}} |k|^{2m - \frac{2}{3}} \| \partial_{y} \theta_{k}\|^2_{L^2}  \Big)^{\frac{1}{2}} \Big( \sum_{k\neq 0}  |k|^{2m } \| \partial_{y} \theta_{k}\|^2_{L^2}  \Big)^{\frac{1}{2}} \Big( \sum_{k\neq 0} |k|\| \psi_{k}\|_{L^\infty}  \Big)\nn\\
			&\lesssim \nu^{-\frac{1}{2}} \mathscr{D}_{\theta, 1}^{\frac{1}{2}} \mathscr{D}_{\omega, 4}^{\frac{1}{2}} \mathcal{E}_{\theta}^{\frac{1}{2}}.
		\end{align}
		For the term $ R_{\neq, \neq, 1} $, we decompose it into the following two parts
		\begin{align}\label{r==1}
			R_{\neq, \neq, 1} &= -2 e^{2\delta_1 \lambda^{\frac{1}{3}} t} \nu^{\frac{1}{3}} \sum_{\substack{ k,\ell \in \Z\setminus\{0\} \\ k \neq \ell }} |k|^{2m - \frac{4}{3}} \operatorname{Re} \Big\langle ik\theta_k,  \partial_{yy} \psi_{k-\ell} \cdot i\ell \theta_\ell  \Big\rangle\nn\\
			&\quad -2 e^{2\delta_1 \lambda^{\frac{1}{3}} t} \nu^{\frac{1}{3}} \sum_{\substack{ k,\ell \in \Z\setminus\{0\} \\ k \neq \ell }} |k|^{2m - \frac{4}{3}} \operatorname{Re} \Big\langle ik\theta_k,  \partial_{y} \psi_{k-\ell} \cdot i\ell \partial_{y}\theta_\ell  \Big\rangle\nn\\
			&:= R_{1,L} + R_{1, R}.
		\end{align}
		For the first term in the above expression, we need to decompose the frequency. Similar to \eqref{p==2}, we have
		\begin{align*}
			R_{1,L} &= -2 e^{2\delta_1 \lambda^{\frac{1}{3}} t} \nu^{\frac{1}{3}} \sum_{\substack{ k,\ell \in \Z\setminus\{0\} \\ k \neq \ell }} \Big( \chi_0(k,\ell) + \chi_1(k,\ell) \Big)  |k|^{2m - \frac{4}{3}} \operatorname{Re} \Big\langle ik\theta_k,  \partial_{yy} \psi_{k-\ell} \cdot i\ell \theta_\ell  \Big\rangle\nn\\
			&:= 	R_{1,L, \chi_0} + R_{1,L, \chi_1}.
		\end{align*}
		For $ (k,\ell)\in A_1$, we have $ |k|\lesssim |\ell|$. Then, applying Young's inequality and the estimate from \eqref{r0=} give rise to
		\begin{align}\label{r1x0}
			R_{1,L, \chi_0} &\lesssim e^{2\delta_1 \lambda^{\frac{1}{3}} t} \nu^{\frac{1}{3}} \sum_{\substack{ k,\ell \in \Z\setminus\{0\} \\ k \neq \ell }} |k|^{m+\frac{1}{3}} \| \theta_{k}\|_{L^2} \|\partial_{yy} \psi_{k-\ell} \|_{L^\infty} |\ell|^{m+\frac{1}{3}} \| \theta_{\ell}\|_{L^2}\nn\\
			&\lesssim \mu^{-\frac{1}{6}} \mathscr{D}_{\theta, 3} \mathcal{E}_{\omega}^{\frac{1}{2}}.
		\end{align}
		For $ (k,\ell)\notin A_1$, we have $ |k|\lesssim |k-\ell|+ |\ell|\lesssim |k-\ell|$.   Using Young's inequality and the fact that $ \|\partial_{yy} \psi_{k}  \|_{L^2} \le \|\omega_{k}  \|_{L^2}$ yields
		\begin{align}\label{r1x1}
			R_{1,L, \chi_1} &\lesssim e^{2\delta_1 \lambda^{\frac{1}{3}} t} \nu^{\frac{1}{3}} \sum_{\substack{ k,\ell \in \Z\setminus\{0\} \\ k \neq \ell }} |k|^{m+\frac{1}{3}} \| \theta_{k}\|_{L^2} |k-\ell|^{m+\frac{1}{3}}\|\partial_{yy} \psi_{k-\ell} \|_{L^2}  \| \theta_{\ell}\|_{L^\infty}\nn\\
			&\lesssim \nu^{\frac{1}{6}} \mu^{-\frac{1}{6}} \mathscr{D}_{\theta, 3}^{\frac{1}{2}} \mathscr{D}_{\omega, 3}^{\frac{1}{2}} \sum_{k\neq 0} \|\theta_{k} \|_{L^\infty} \nn\\
			&\lesssim \nu^{\frac{1}{6}} \mu^{-\frac{1}{6}} \mathscr{D}_{\theta, 3}^{\frac{1}{2}} \mathscr{D}_{\omega, 3}^{\frac{1}{2}} \sum_{k\neq 0} |k|^{\frac{1}{6}-m} \cdot |k|^{\frac{m}{2}}\|\theta_{k} \|^{\frac{1}{2}}_{L^2} \cdot |k|^{\frac{m}{2}-\frac{1}{6}}\|\partial_{y}\theta_{k} \|^{\frac{1}{2}}_{L^2} \nn\\
			&\lesssim \mu^{-\frac{1}{6}}\mathscr{D}_{\theta, 3}^{\frac{1}{2}} \mathscr{D}_{\omega, 3}^{\frac{1}{2}} \mathcal{E}_{\theta}^{\frac{1}{2}},
		\end{align}
		where we have used the fact that $ m> \frac{2}{3}$. Next, we handle the most challenging term $ R_{1,R}$. First, we need to decompose the frequency $k$, obtaining
		\begin{align*}
			R_{1,R} &= -2 e^{2\delta_1 \lambda^{\frac{1}{3}} t} \nu^{\frac{1}{3}} \sum_{\substack{ k,\ell \in \Z\setminus\{0\} \\ k \neq \ell }} |k|^{m - \frac{2}{3}}\Big[|k|^{m - \frac{2}{3}} - |\ell|^{m - \frac{2}{3}} \Big] \operatorname{Re} \Big\langle ik\theta_k,  \partial_{y} \psi_{k-\ell} \cdot i\ell \partial_{y}\theta_\ell  \Big\rangle\nn\\
			&\quad -2 e^{2\delta_1 \lambda^{\frac{1}{3}} t} \nu^{\frac{1}{3}} \sum_{\substack{ k,\ell \in \Z\setminus\{0\} \\ k \neq \ell }} |k|^{m - \frac{2}{3}}\cdot |\ell|^{m - \frac{2}{3}}  \operatorname{Re} \Big\langle ik\theta_k,  \partial_{y} \psi_{k-\ell} \cdot i\ell \partial_{y}\theta_\ell  \Big\rangle := R_{1,R, L} + R_{1,R,R}.
		\end{align*}
		For the second term in the above expression, we directly obtain the result using Young's inequality and \eqref{r0=}:
		\begin{align}\label{r1rr}
			 R_{1,R,R} &= - e^{2\delta_1 \lambda^{\frac{1}{3}} t} \nu^{\frac{1}{3}} \sum_{\substack{ k,\ell \in \Z\setminus\{0\} \\ k \neq \ell }} |k|^{m - \frac{2}{3}}\cdot |\ell|^{m - \frac{2}{3}}  \operatorname{Re} \Big\langle ik\theta_k,  \partial_{yy} \psi_{k-\ell} \cdot i\ell \theta_\ell  \Big\rangle\nn\\
			 &\lesssim \mathscr{D}_{\theta, 3} \sum_{k\neq 0} \| \partial_{yy} \psi_{k}\|_{L^\infty} \lesssim \mu^{-\frac{1}{6}}\mathscr{D}_{\theta, 3} \mathcal{E}_{\omega}^{\frac{1}{2}}.
		\end{align}
		For the first term, by applying a method similar to that used in \eqref{p==2}, we obtain
		\begin{align*}
			R_{1,R, L} &= -2 e^{2\delta_1 \lambda^{\frac{1}{3}} t} \nu^{\frac{1}{3}} \sum_{\substack{ k,\ell \in \Z\setminus\{0\} \\ k \neq \ell }} \Big( \chi_0(k,\ell) + \chi_1(k,\ell) \Big) |k|^{m - \frac{2}{3}}\nn\\
&\qquad\times\Big[|k|^{m - \frac{2}{3}} - |\ell|^{m - \frac{2}{3}} \Big] \operatorname{Re} \Big\langle ik\theta_k,  \partial_{y} \psi_{k-\ell} \cdot i\ell \partial_{y}\theta_\ell  \Big\rangle\nn\\
			& := R_{1,R, \chi_0} + R_{1,R, \chi_1}.
		\end{align*}
		We first consider the first term in the above expression.  For $ (k,\ell)\in A_1$, we have $ |k|\approx |\ell|$. By applying the mean value theorem, we obtain
		\begin{align*}
			\Big| |k|^{m - \frac{2}{3}} - |\ell|^{m - \frac{2}{3}}  \Big| = (m-\frac{2}{3}) |k-\ell| \Big| \alpha_1 k + (1-\alpha_1) \ell\Big|^{m-\frac{2}{3}-1} \lesssim |k-\ell| | \ell|^{m-\frac{5}{3}}.
		\end{align*}
		Applying Young's inequality, the Gagliardo-Nirenberg inequality, and the fact from  \eqref{w_fact} give rise to
		\begin{align}\label{r1rx0}
			R_{1,R, \chi_0} &\lesssim e^{2\delta_1 \lambda^{\frac{1}{3}} t} \nu^{\frac{1}{3}} \sum_{\substack{ k,\ell \in \Z\setminus\{0\} \\ k \neq \ell }} |k|^{m+\frac{1}{3}} \|\theta_{k} \|_{L^2} \cdot |k-\ell|^{\frac{2}{3}} \| \partial_{y} \psi_{k-\ell}\|_{L^\infty} \cdot |\ell|^{m-\frac{1}{3}} \| \partial_{y} \theta_{\ell}\|_{L^2}\nn\\
			&\lesssim e^{2\delta_1 \lambda^{\frac{1}{3}} t} \nu^{-\frac{1}{6}} \Big( \sum_{k\neq 0} \nu^{\frac{1}{3}} |k|^{2m+\frac{2}{3}} \|\theta_{k} \|_{L^2}^2\Big)^{\frac{1}{2}} \Big( \sum_{k\neq 0} \nu^{\frac{2}{3}} |k|^{2m-\frac{2}{3}} \|\partial_{y}\theta_{k} \|_{L^2}^2\Big)^{\frac{1}{2}} \sum_{k\neq 0} |k|^{\frac{2}{3}} \| \partial_{y} \psi_{k}\|_{L^\infty}\nn\\
			&\lesssim \nu^{-\frac{1}{6}} \mathscr{D}_{\theta, 3}^{\frac{1}{2}} \mathcal{E}_{\theta}^{\frac{1}{2}} \sum_{k\neq 0} |k|^{\frac{1}{6}} \| \omega_k\|_{L^2}\nn\\
			&\lesssim \nu^{-\frac{1}{6}} \mu^{-\frac{1}{6}} \mathscr{D}_{\theta, 3}^{\frac{1}{2}} \mathscr{D}_{\omega, 3}^{\frac{1}{2}} \mathcal{E}_{\theta}^{\frac{1}{2}}.
		\end{align}
		Similarity, we have
		\begin{align}\label{r1rx1}
			R_{1,R, \chi_1} &\lesssim e^{2\delta_1 \lambda^{\frac{1}{3}} t} \nu^{\frac{1}{3}} \sum_{\substack{ k,\ell \in \Z\setminus\{0\} \\ k \neq \ell }} |k|^{m+\frac{1}{3}} \|\theta_{k} \|_{L^2} \cdot |k-\ell|^{m+\frac{1}{3}} \| \partial_{y} \psi_{k-\ell}\|_{L^\infty} \cdot \| \partial_{y} \theta_{\ell}\|_{L^2}\nn\\
			&\lesssim e^{2\delta_1 \lambda^{\frac{1}{3}} t} \nu^{ -\frac{1}{3}} \Big( \sum_{k\neq 0} \nu^{\frac{1}{3}} |k|^{2m+\frac{2}{3}} \|\theta_{k} \|_{L^2}^2\Big)^{\frac{1}{2}} \Big( \sum_{k\neq 0}  |k|^{2m+\frac{2}{3}} \|\partial_{y}\psi_{k} \|_{L^\infty}^2\Big)^{\frac{1}{2}} \sum_{k\neq 0} \nu^{ \frac{1}{2}}  \| \partial_{y} \theta_{k}\|_{L^2}\nn\\
			&\lesssim  \nu^{- \frac{1}{3}} \mathscr{D}_{\theta, 3}^{\frac{1}{2}} \mathscr{D}_{\theta, 1}^{\frac{1}{2}} \mathcal{E}_{\omega}^{\frac{1}{2}},
		\end{align}
		where we have used the  fact that
		\begin{align*}
			e^{2\delta_1 \lambda^{\frac{1}{3}} t} \sum_{k\neq 0}  |k|^{2m+\frac{2}{3}} \|\partial_{y}\psi_{k} \|_{L^\infty}^2 &\lesssim e^{2\delta_1 \lambda^{\frac{1}{3}} t}	\sum_{k\neq 0}  |k|^{2m+\frac{2}{3}} \|\partial_{y}\psi_{k} \|_{L^2}  \|\partial_{yy}\psi_{k} \|_{L^2} \\
			&\lesssim e^{2\delta_1 \lambda^{\frac{1}{3}} t} \sum_{k\neq 0}  |k|^{2m-\frac{1}{3}}  \| \omega_k\|^2_{L^2} \lesssim \mathcal{E}_{\omega}.
		\end{align*}
		Collecting \eqref{r0=}-\eqref{r1rx1}, we can complete the proof of the lemma.
	\end{proof}
\end{lemma}

\begin{lemma}\label{le_Ewta}
	Under the conditions of Theorem \ref{thm1.2}, then for any integer $k \neq 0$ and any $t>0$, we have
	\begin{align}
		\frac{d}{dt} \mathcal{E}_{\theta} + \frac{1}{4} \mathscr{D}_{\theta} &\lesssim  \nu^{-\frac{1}{2}} \mathscr{D}_{\theta}^{\frac{1}{2}} \mathscr{D}_{\omega}^{\frac{1}{2}} \mathcal{E}_{\theta}^{\frac{1}{2}} + (\nu^{-\frac{1}{2}}  + \mu^{-\frac{1}{6}}) \mathscr{D}_{\theta}  \mathcal{E}_{\omega}^{\frac{1}{2}} + \nu^{-\frac{1}{6}} \mu^{-\frac{1}{3}} \mathscr{D}_{\theta}^{\frac{1}{2}} \mathscr{D}_{\omega}^{\frac{1}{2}} \mathcal{E}_{\theta}^{\frac{1}{2}},\label{Eta}\\
		\frac{d}{dt} \mathcal{E}_{\omega} + \frac{1}{16} c_{\beta} \mathscr{D}_{\omega} &\lesssim \mu^{-\frac{1}{2}}  \mathscr{D}_{\omega}  \mathcal{E}_{\omega}^{\frac{1}{2}} + \mu^{-\frac{1}{3}} \nu^{-\frac{2}{3}} \mathscr{D}_{\theta}.\label{Ew}
	\end{align}
	\begin{proof}
	First, we take $ \delta_1 = \min\{ \frac{1}{8}, \frac{c_0}{4}\}$ (where $c_0$ is a constant determined in Lemma \ref{le_w7}). Then,  for $0 < \mu, \nu \le 1$, we can directly obtain \eqref{Eta} using Lemmas \ref{le_w7}-\ref{le_nonta3}. Next, we proceed to prove \eqref{Ew}. It suffices to estimate the nonlinear terms $ \mathcal{N}_{\omega, \alpha}$, $ \mathcal{N}_{\omega, \beta}$, and $ \mathcal{N}_{\omega, \gamma}$. For the term $ \mathcal{N}_{\omega, \alpha}$	, using a method similar to Lemma \ref{le_nonta1} and the 	boundedness of the operators $\mathfrak{J}_k$ and $ [\partial_{y}, \mathfrak{J}_k]$ (from Lemma \ref{le2.1} and Lemma \ref{le2.2}), we obtain
		\begin{align}\label{nonw1}
			\mathcal{N}_{\omega, \alpha} \lesssim \mu^{-\frac{1}{2}} \mathscr{D}_{\omega}  \mathcal{E}_{\omega}^{\frac{1}{2}}.
		\end{align}
		For the term $  \mathcal{N}_{\omega, \beta}$, using a method similar to Lemma \ref{le_nonta2} and the fact that
 $$ \|\partial_{y} \mathfrak{J}_k f \|_{L^2} \lesssim \|[ \partial_{y} , \mathfrak{J}_k] f \|_{L^2} + \| \mathfrak{J}_k \partial_{y}f \|_{L^2} \lesssim \| \nabla_k f\|_{L^2},$$ we obtain
		\begin{align}\label{nonw2}
				\mathcal{N}_{\omega, \beta} \lesssim \mu^{-\frac{1}{2}} \mathscr{D}_{\omega}  \mathcal{E}_{\omega}^{\frac{1}{2}}.
		\end{align}
		For the term $ \mathcal{N}_{\omega, \gamma}$, we can directly obtain the result by following the method in Lemma \ref{le_nonta3}:
		\begin{align}\label{nonw3}
		\mathcal{N}_{\omega, \gamma} \lesssim \mu^{-\frac{1}{2}} \mathscr{D}_{\omega}  \mathcal{E}_{\omega}^{\frac{1}{2}}.
		\end{align}
		Combining the above three equations with Lemma \ref{le_w7} and \eqref{lw=}, we obtain \eqref{Ew}.
	\end{proof}
\end{lemma}

\subsection{Completion the proof of the nonlinear stability of Theorem \ref{thm1.2}}
In this subsection, we investigate the enhanced dissipation and the inviscid damping of the solution . Before proceeding with the proof, for convenience, we define the following two energy functionals
\begin{align*}
	\mathcal{E}_{\theta, \text{total}}(t) := \mathcal{E}_{\theta}(t) + \frac{1}{4}\int_{0}^{t} \mathscr{D}_{\theta}(\tau) \mathrm{d} \tau, \quad \mathcal{E}_{\omega, \text{total}}(t) := \mathcal{E}_{\omega}(t) + \frac{1}{16}c_{\beta}\int_{0}^{t} \mathscr{D}_{\omega}(\tau) \mathrm{d} \tau.
\end{align*}
Then, it follows from Lemma \ref{le_Ewta} that
\begin{align}
	 &\mathcal{E}_{\theta, \text{total}}(t) \lesssim \mathcal{E}_{\theta, \text{total}}(0) + \big(\nu^{-\frac{1}{2}} + \nu^{-\frac{1}{6}}\mu^{-\frac{1}{3}} \big) \mathcal{E}_{\theta, \text{total}}(t) \mathcal{E}_{\omega, \text{total}}^{\frac{1}{2}}(t),\label{Eta_total}\\
	 &\mathcal{E}_{\omega, \text{total}}(t) \lesssim \mathcal{E}_{\omega, \text{total}}(0) +  \mu^{-\frac{1}{2}} \mathcal{E}_{\omega, \text{total}}^{\frac{3}{2}}(t) + \mu^{-\frac{1}{3}} \nu^{-\frac{2}{3}} 	\mathcal{E}_{\theta, \text{total}}(t).\label{Ew_total}
\end{align}
Moreover, it is straightforward to verify from
$$  \sum_{0\le j\le 1 } \| (\mu^{\frac{1}{3} } \partial_{y})^j \langle \partial_x\rangle^{m-\frac{j}{3}} \omega_{in} \|_{L^2} \le \varepsilon_0 \min\big\{ \mu^{-\frac{1}{2} }, \nu^{-\frac{1}{2} }\big\},$$ and  $$  \sum_{0\le j\le 1 } \| (\nu^{\frac{1}{3} } \partial_{y})^j \langle \partial_x\rangle^{m-\frac{j}{3}} \theta_{in} \|_{L^2} \le \varepsilon_1 \min\big\{ \mu^{-1 }, \nu^{-1} \big\},$$
 that $ \mathcal{E}_{\theta, \text{total}}(0) \le C \varepsilon_1 \min\big\{ \mu^{-2 }, \nu^{-2} \big\}$ and $ \mathcal{E}_{\omega, \text{total}}(0) \le C \varepsilon_0 \min\big\{ \mu^{-1 }, \nu^{-1} \big\}$.

Therefore, for sufficiently small $ \varepsilon_0$ and $\varepsilon_1 $, by the bootstrap argument, we can deduce from \eqref{Eta_total} and \eqref{Ew_total} that
\begin{align*}
	\mathcal{E}_{\theta, \text{total}}(t) \le C \varepsilon_1 \min\big\{ \mu^{-2 }, \nu^{-2} \big\}, \quad \mathcal{E}_{\omega, \text{total}}(t) \le C \varepsilon_0 \min\big\{ \mu^{-1 }, \nu^{-1} \big\}.
\end{align*}

This completes the proof of Theorem \ref{thm1.2}. \hspace{7.8cm}
$\square$

\bigskip

\section*{Appendix}

\setcounter{theorem}{0}
\setcounter{equation}{0}
\renewcommand{\theremark}{A.\arabic{remark}}
\renewcommand{\theequation}{A.\arabic{equation}}
\renewcommand{\thetheorem}{A.\arabic{theorem}}

In this section, we provide some useful properties of integral kernel functions. Inspired by  the work in \cite{WangCMP2025},  we define the integral form $h(y,y')$ as follows $$ h(y, y'):= \frac{1}{y-y'} \int_{y}^{y'} U'''(s) (y' - s)^2 \mathrm{d} s.$$
\begin{lemma}\label{A1}
	Consider $W(t,y)$ satisfying the initial and boundary conditions $W_{in}|_{y = \pm 1} = 0$. If the initial condition satisfies \eqref{W0}, then the function $ h(y,y')$ has the following estimate
	\begin{align}\label{1h}
		 |h(y,y')| \lesssim \delta_0 |y-y'|^2, \ \  |\partial_{y} h(y, y')| + |\partial_{y'} h(y, y')| \lesssim \delta_0 |y-y'|,\ \  |\partial_{y} \partial_{y'} h(y, y')| \lesssim \delta_0.
	\end{align}
	\begin{proof}
		First, by applying the Mean Value Theorem for integrals, we directly obtain
		\begin{align}\label{h1}
			h(y,y') &= -\alpha^2 U'''(\alpha y + (1-\alpha)y') (y-y')^2 \nn\\
			&\le \| \partial_{yy} W\|_{L^\infty} (y-y')^2 \lesssim \delta_0 |y-y'|^2, \quad  \alpha \in (0,1).
		\end{align}
		Then, by directly computing the spatial derivative, we can get
		\begin{align}\label{h2}
			\partial_{y} h(y, y') = -U'''(y)(y-y') + \alpha^2 (y-y') U'''( \alpha y + (1-\alpha)y')  \lesssim \delta_0 |y-y'|, \quad \quad \alpha \in (0,1).
		\end{align}
		Similarly, we have $ \partial_{y'} h(y, y') \lesssim \delta_0 |y-y'|$ and $ \partial_{y'} \partial_{y} h(y, y') \lesssim \delta_0$.
	\end{proof}
\end{lemma}

Motivated by \cite{WangCMP2025},  we introduce the following two integral kernel functions $ \mathscr{K}_1(y,y') = (y-y')U''(y) G_k(y,y')$ and  $\mathscr{K}_2 (y,y') = h(y,y') G_k(y,y')$, whose properties will be rigorously proved in the subsequent lemma.

\begin{lemma}\label{A2}
	Under the conditions of Lemma \ref{A1}, for any $ k \neq 0$, we have the following
	\begin{align}\label{1k}
		\| \mathscr{K}_1(y, y')\|_{L^2_{y,y'}} \lesssim  \frac{\delta_0}{k^2} ,
	\end{align}
	and
	\begin{align}\label{2k}
		\| \partial_{y}\mathscr{K}_1(y, y')\|_{L^2_{y,y'}} + \| \partial_{y'}\mathscr{K}_1(y, y')\|_{L^2_{y,y'}} \lesssim \frac{\delta_0}{|k|} ,
	\end{align}
	and
\begin{align}\label{3k}
\| \partial_{y} \partial_{y'} \mathscr{K}_1(y, y')\|_{L^2_{y,y'}}  \lesssim  \delta_0.
\end{align}
	\begin{proof}
		First, by directly computing, we obtain
		\begin{align}\label{k1}
			\| \mathscr{K}_1\|^2 _{L^2_{y,y'}} &= \int_{-1}^{1} \int_{-1}^{1} |U''(y)|^2 |y-y'|^2 |G_{k}(y, y')|^2  \mathrm{d} y' \mathrm{d} y\nn\\
			&\le \| U''\|^2_{L^\infty} \Big( \underbrace{ \int_{-1}^{1} \int_{-1}^{y}  |y-y'|^2 \frac{1}{k^2 \sinh^2 (2 k)} \big[ \sinh (k(1-y)) \sinh \left(k\left(1+y^{\prime}\right)\right)\big]^2 \mathrm{d} y' \mathrm{d} y}_{I_1} \nn\\
			&\quad + \underbrace{ \int_{-1}^{1} \int_{y}^{1}  |y-y'|^2 \frac{1}{k^2 \sinh^2 (2 k)} \big[ \sinh (k(1+y)) \sinh \left(k\left(1-y^{\prime}\right)\right)\big]^2 \mathrm{d} y' \mathrm{d} y}_{I_2} \Big).
		\end{align}
		Next, we estimate each term of the two expressions above term by term,
		\begin{align}\label{I1}
			I_1 & = \frac{1}{k^2 \sinh^2 (2k)} \int_{-1}^1 [\sinh k(1-y)]^2 \int_{-1}^y \left(y-y^{\prime}\right)^2 \left[-1+\frac{e^{2 k\left(1+y^{\prime}\right)}+e^{-2k\left(1+y^{\prime}\right)}}{4}\right] \mathrm{d} y^{\prime} \mathrm{d} y\nn\\
			& = \frac{1}{k^2 \sinh^2 (2k)} \int_{-1}^1 [\sinh k(1-y)]^2   \left[\left. \frac{1}{3} \left(y^{\prime}-y\right)^3\right|_{-1} ^y + \int_{-1}^y (y-y')^2 \frac{\cosh 2k(1+y')}{2} \mathrm{d} y^{\prime} \right]  \mathrm{d} y\nn\\
			& := I_{11} + I_{12}.
		\end{align}
		For the term $ I_{11}$, we decompose it into the following two terms:
		\begin{align}\label{I11}
			I_{11} = \frac{1}{3 k^2 \sinh ^2 2k} \int_{-1}^1\left[-1+\frac{\cosh 2k(1-y)}{2}\right](y+1)^3 \mathrm{d} y := I_{111} + I_{112}.
		\end{align}
			Term $ I_{111}$ is obtained by directly applying integration by parts, yielding
		\begin{align}\label{I111}
				\big| I_{111}  \big| \lesssim \big| \frac{4}{3k^2 \sinh ^2 2k}\big| \lesssim \frac{1}{k^4}.
		\end{align}
		Similarly, term $I_{112}$ can also be directly computed, yielding
		\begin{align}\label{I112}
			\big| I_{112}  \big| \lesssim \big| \frac{1}{k^4 \sinh ^2 2k} + \frac{\cosh 4k}{k^4 \sinh^2 2k} \big| \lesssim \frac{1}{k^4}.
		\end{align}
		By combining \eqref{I111} and \eqref{I112}, we obtain \begin{align}\label{I11_es}
			|I_{11}| \lesssim \frac{1}{k^4}.
		\end{align}
		Next, we estimate term $I_{12}$. First, we consider the integral terms with respect to $y'$ in $ I_{12}$, which are given by
		\begin{align}\label{I12_re}
			\int_{-1}^y (y-y')^2 \frac{\cosh 2k(1+y')}{2} \mathrm{d} y^{\prime} & = -\frac{1}{k} \int_{-1}^y\left(y^{\prime}-y\right) \cdot \sinh  2k\left(1+y^{\prime}\right) d y^{\prime} \nn\\
			& = - \frac{1}{2k^2} (y + 1) + \frac{1}{4k^3} \sinh 2k(1+y).
		\end{align}
		Therefore, we can decompose term $I_{12}$ into the following two terms
		\begin{align}\label{AI122}
			I_{12} &= \frac{1}{2k^2 \sinh ^2 2k} \int_{-1}^1 [\sinh k(1-y)]^2   \left[ - \frac{1}{2k^2} (y + 1) + \frac{1}{4k^3} \sinh 2k(1+y) \right]  \mathrm{d} y \nn\\
			& := I_{121} + I_{122}.
		\end{align}
		For the terms $ I_{121}$ and $I_{122}$ , we can directly compute it to obtain
		\begin{align*}
			\big| I_{121}\big| \lesssim  \frac{1}{4k^4\sinh ^2 2k} \int_{-1}^{1} 1 + \cosh 2k(1-y)  \mathrm{d} y \lesssim \frac{1}{k^4},
		\end{align*}
		and
		\begin{align*}
			I_{122} &= -\frac{1}{8k^5 \sinh^2 2k} \int_{-1}^{1} \sinh 2k(1+y) \mathrm{d} y + \frac{1}{16k^5 \sinh^2 2k} \int_{-1}^{1} \cosh 2k(1-y) \sinh 2k(1+y) \mathrm{d} y \nn\\
			& \lesssim \frac{1}{k^4} + \frac{1}{16k^5 \sinh^2 2k} \int_{-1}^{1} \frac{\sinh 4k}{2} + \frac{\sinh 4ky}{2} \mathrm{d} y  \lesssim \frac{1}{k^4}.
		\end{align*}
		By combining the above two equations and \eqref{I11_es}, we obtain
		\begin{align}\label{I1_es}
		|I_{1}| \lesssim \frac{1}{k^4}.
		\end{align}
		For term $I_2$, we obtain the result by handling it in a similar manner to $I_1$, yielding
		\begin{align}\label{I2_es}
			|I_{2}| \lesssim \frac{1}{k^4}.
		\end{align}
		By substituting \eqref{I1_es} and \eqref{I2_es} into \eqref{k1}, we obtain \eqref{1k}.
		Next, we estimate terms $ \partial_{y} \mathscr{K}_1$,  $ \partial_{y'} \mathscr{K}_1$ and $ \partial_{y} \partial_{y'} \mathscr{K}_1$. By directly differentiating, we obtain
		\begin{align}
			\partial_{y} \mathscr{K}_1 &= U'''(y) (y-y') G_k(y,y') + U''(y) G_k(y,y') + U''(y) (y-y') \partial_{y}G_k(y,y'),\label{diffy}\\
			\partial_{y'} \mathscr{K}_1 &=- U''(y) G_k(y,y') + U''(y) (y-y') \partial_{y'}G_k(y,y'),\label{diffy'}\\
			 \partial_{y} \partial_{y'} \mathscr{K}_1 &= - U'''(y) G_k(y,y') - U''(y)  \partial_{y}G_k(y,y') + U'''(y) (y-y') \partial_{y'}G_k(y,y')\nn\\
			&\quad + U''(y)  \partial_{y'}G_k(y,y') + U''(y) (y-y') \partial_{yy'}G_k(y,y').
		\end{align}
		For the first term on the right-hand side of \eqref{diffy}, we can directly obtain it from the proof of $ \| \mathscr{K}_1\|^2 _{L^2_{y,y'}}$. For the second term on the right-hand side of \eqref{diffy}, integration by parts introduces an extra factor of $\frac{1}{k}$. However, for the third term, since $ \partial_{y} G_k(y,y')$ cancels out a factor of $\frac{1}{k}$, the remaining estimates follow from those in $ \| \mathscr{K}_1\|^2 _{L^2_{y,y'}}$. Therefore, we obtain the following estimates  (a similar estimate holds for terms $ \partial_{y'} \mathscr{K}_1 $ and $ \partial_{y}\partial_{y'} \mathscr{K}_1$)
		\begin{align*}
			&\| \partial_{y}\mathscr{K}_1(y, y')\|_{L^2_{y,y'}} + \| \partial_{y'}\mathscr{K}_1(y, y')\|_{L^2_{y,y'}} \lesssim \|U''' \|_{L^\infty} \frac{1}{k^2} + \|U'' \|_{L^\infty} \frac{1}{|k|},\\
			& \| \partial_{y} \partial_{y'}\mathscr{K}_1(y, y')\|_{L^2_{y,y'}} \lesssim \|U''' \|_{L^\infty} \frac{1}{k^2} + \|(U'', U''') \|_{L^\infty} \frac{1}{|k|} + \|U'' \|_{L^\infty},
		\end{align*}
		which uses the fact that $ \|(U'', U''') \|_{L^\infty} \lesssim \| W(t,\cdot)\|_{H^4_y} \le \delta_0$, completing the proof of Lemma \ref{A2}.
	\end{proof}
\end{lemma}
Additionally, for the estimate of the integral kernel $ \mathscr{K}_2(y, y')$, using the results of Lemma \ref{A1} and Lemma \ref{A2}, we can directly obtain the following result.
\begin{corollary}\label{A3}
	Under the conditions of Lemma \ref{A1}, for any $ k \neq 0$, we have the following
	\begin{align*}
	\| \mathscr{K}_2(y, y')\|_{L^2_{y,y'}} &\lesssim  \frac{\delta_0}{k^2} ,\\
	\| \partial_{y}\mathscr{K}_2(y, y')\|_{L^2_{y,y'}} + \| \partial_{y'}\mathscr{K}_2(y, y')\|_{L^2_{y,y'}} &\lesssim \frac{\delta_0}{|k|} ,\\
	\| \partial_{y} \partial_{y'} \mathscr{K}_2(y, y')\|_{L^2_{y,y'}}  &\lesssim  \delta_0.
	\end{align*}
\end{corollary}

 \section*{Acknowledgement} {This work is  supported by the Guangdong Provincial Natural Science Foundation under grant  2024A1515030115.}

 \vskip .1in
\noindent{\bf Data Availability Statement} Data sharing is not applicable to this article as no
data sets were generated or analysed during the current study.

\vskip .1in

\noindent{\bf Conflict of Interest} The authors declare that they have no conflict of interest. The
authors also declare that this manuscript has not been previously published, and
will not be submitted elsewhere before your decision.

\end{document}